\documentclass[a4paper,fleqn]{cas-sc}

\usepackage[authoryear,longnamesfirst]{natbib}

\def\tsc#1{\csdef{#1}{\textsc{\lowercase{#1}}\xspace}}
\tsc{WGM}
\tsc{QE}

\usepackage{lipsum}
\usepackage{amsfonts}
\usepackage{bm}   
\usepackage{amsmath}
\usepackage{graphicx}
\usepackage{dsfont}
\usepackage{epstopdf}
\usepackage{algorithm}
\usepackage{algpseudocode}
\usepackage{booktabs}
\usepackage{enumerate}
\usepackage[toc]{appendix}
\usepackage{multirow}
\usepackage{soul}
\usepackage{caption}
\usepackage{subcaption}
\usepackage{scalerel,stackengine}
\usepackage{tikz}
\usepackage{mathtools}
\usepackage{amsopn}
\usepackage{algpseudocode}
\tikzset{
  treenode/.style = {shape=rectangle, rounded corners,
                     draw, align=center,
                     top color=white, bottom color=blue!20},
  root/.style     = {treenode, font=\Large, bottom color=red!30},
  env/.style      = {treenode, font=\ttfamily\normalsize},
  dummy/.style    = {circle,draw}
}

\newcommand{\reallywidehat}[1]{%
\savestack{\tmpbox}{\stretchto{%
\scaleto{%
\scalerel*[\widthof{\ensuremath{#1}}]                        
{\kern-.6pt\bigwedge\kern-.6pt}%
{\rule[-\textheight/2]{1ex}{\textheight}}
}{\textheight}%
}{0.5ex}}%
\ensurestackMath{\stackon[1pt]{#1}{\tmpbox}}%
}

\ifpdf
  \DeclareGraphicsExtensions{.eps,.pdf,.png,.jpg}
\else
  \DeclareGraphicsExtensions{.eps}
\fi

\newtheorem{prop}{Proposition}

\def\RR{\mathbb{R}}

\DeclareMathOperator*{\argmax}{arg\,max}

\DeclareMathOperator*{\argmin}{arg\,min}

\DeclarePairedDelimiter\floor{\lfloor}{\rfloor}
\newcommand{\ep}{\varepsilon_{\mathcal{T}}}
\newcommand{\R}{\mathbb R}

\definecolor{greenf}{rgb}{.054, .5, .005}
\begin{document}
\let\WriteBookmarks\relax
\def\floatpagepagefraction{1}
\def\textpagefraction{.001}

\shorttitle{Approximation of Optimal Control Problems for the NS equation via multilinear HJB-POD}    

\shortauthors{M. Falcone, G. Kirsten and L. Saluzzi}  

\title [mode = title]{Approximation of Optimal Control Problems for the Navier-Stokes equation via multilinear HJB-POD}  



%

\author[1]{Maurizio Falcone}[type=editor]


\affiliation[1]{email={falcone@mat.uniroma1.it,}, organization={Universit\`a di Roma La Sapienza},
            addressline={P. Aldo Moro, 5}, 
            city={Roma},
            postcode={00185}, 
            state={},
            country={Italy}}

\author[2]{Gerhard Kirsten}[type=editor]

%

%
%
%
\affiliation[2]{email={g.kirsten@unibo.it,}, organization={Dipartimento di Matematica, Universit\`a di Bologna },
            addressline={Piazza P.ta S. Donato}, 
            city={Bologna},
            postcode={51000}, 
            country={Italy}}
%
\cortext[1]{Corresponding author}
%


\author[3]{Luca Saluzzi}[type=editor]

\cormark[1]

%
%
\affiliation[3]{email={l.saluzzi@ic.ac.uk,},
organization={Department of Mathematics, Imperial College London},
            addressline={South Kensington Campus, SW7 2AZ}, 
           city={London},
            country={United Kingdom}}
%
%


\begin{abstract}
We consider the approximation of some optimal control problems for the Navier-Stokes equation via a Dynamic Programming approach. These control problems arise in many industrial applications and
are very challenging from the numerical point of view since the semi-discretization of the dynamics corresponds to an evolutive system of ordinary differential equations in very high dimension. The typical approach is based on the Pontryagin maximum principle and leads to a two point boundary value problem. Here we present a different approach based on the value function and the solution of a Bellman, a challenging problem in high dimension. We mitigate the curse of dimensionality via a recent multilinear approximation of the dynamics coupled with a dynamic programming scheme on a tree structure. We discuss several aspects related to the implementation of this new approach and we present some numerical examples to illustrate the results on classical control problems studied in the literature.
\end{abstract}



\begin{keywords}
 \sep dynamic programming \sep optimal control \sep tree structure \sep model order reduction 
\end{keywords}

\maketitle

\section{Introduction}\label{sec:intro}
The control of fluids is an important issue in many industrial problems, e.g. in aerospace and naval industries. The approximation of the fluid around complex geometries usually requires a very careful construction of the grid and is based on finite elements or finite differences/volumes schemes (see e.g. \cite{pir89}, \cite{strikwerda2004} and the references therein). This is known to be a huge computational problem so that model reduction techniques are often applied to
compute the solution \cite{BCOW.17, benner2015review} and study the physical properties of the flow varying various parameters
(e.g. bifurcation phenomena as in \cite{QR2007,SR2018,pichi2022driving}).
A typical example is given by the Navier-Stokes equation for incompressible fluids that we use here as our model problem.\\
 These problems have been studied by many authors from the theoretical point of view analyzing  controllability properties of the system via Carlemann estimates, the interested reader can find in   \cite{fi96, puel2014} a comprehesive presentation of these results. \\
 For optimal control problems the numerical approximation of the Navier-Stokes equation is the starting point. In fact we want to solve a  huge optimization problem where the controlled solution should minimize a cost functional, e.g. to minimize the drag or to stay close to some reference solution. Many authors have contributed to these problems following  the pioneering work of J.L. Lions in \cite{li68} and the numerical approach is mainly based on  the Pontryagin Maximum Principle (PMP) that gives a necessary condition characterizing the optimal couple trajectory/control of the problem.  The numerical solution of the PMP leads to a two point boundary value problem that is feasible in high-dimension but typically produces open-loop controls (see \cite{BCD97}).  In practice this approach can be difficult to implement since it requires a starting guess for the optimal trajectory and for the optimal control (the co-state) that it not available, in particular the co-state is usually hard to initialize. The interested reader can find in in the book  \cite{tr2010} and in the lecture notes \cite{casas06} a general presentation of the results. We just recall that for the PMP approach two different strategies have been proposed: "optimize then discretize" and "discretize then optimize". The first is based on the discretization of the system of optimality conditions obtained for the continuous problem whereas the second starts with the discretization of the optimal control problem and then solves the optimality condition for the finite dimensional problem (see \cite{HPUU2009} for a general presentation of these numerical strategies). \\
As we said, here we follow the Dynamic Programming (DP) approach based on the characterization of the value function as the unique solution of a Hamilton-Jacobi-Bellman equation. This approach is more interesting since it produces a characterization of  optimal controls in feedback form via the knowledge of the value function, however its application to the control of PDEs  has been very limited  due to the "curse of dimensionality". In fact,  adopting the strategy "discretize then optimize", we need to solve the Hamilton-Jacobi-Bellman equation in high-dimension (i.e. the dimension of the discrete state space after the semi-discretization of the continuous problem). It is known that that  nonlinear partial differential equation gives the characterization of the value function as its unique viscosity 
solution in many optimal control problems (see e.g. \cite{BCD97}).  It is interesting to note that this problem is difficult  also in low dimension since the value function is only Lipschitz continuous also when the dynamics and the cost are assumed to be very regular, but in low dimension several methods have been proposed ranging from finite difference methods \cite{Sethian1999}, semi-lagrangian schemes \cite{FF14} and finite volumes. Here we propose a new method for the numerical solution of control problems of Navier-Stokes equations based on the DP approach. The novelty is in the technique used to mitigate the curse of dimensionality via the  coupling of two recent methods: a multilinear approximation of the NS  equation developed in \cite{Kirsten.Simoncini.arxiv2020} and the dynamic programming method for the finite horizon problem on a tree originally developed for nonlinear ordinary differential equations in \cite{afs19} obtaining also a-priori error estimates in \cite{saf21}.
The first method allows to produce a numerical solution in a very compact form via tensor notations whereas the Tree Structure Algorithm (TSA) exploits the compact representation of the systems  and can be coupled with a model reduction approach based on Proper Orthogonal Decomposition (POD).
In fact the tree structure method is rather flexible, we refer to \cite{alla2020} and  \cite{afs20} for recent developments including high-order approximation, the coupling with model reduction techniques, problems with state constraints. 
To set this paper into perspective, let us also mention that the coupling the HJB equation with POD for the approximation of optimal control problems with PDE constraints has been proposed by Kunisch and co-authors in a series of papers \cite{KVX2004, KX2005,kunisch2010optimal} (see also  \cite{HV2005}). They have analyzed optimal control problems mainly for linear parabolic equation and the Burgers equation. The numerical method proposed here is different since it is based on the above mentioned building blocks allowing to mitigate the "curse of dimensionality".
Other techniques have been introduced in the last decades in this direction, among them we mention in particular  sparse grids \cite{GK16} and  tensor decomposition techniques \cite{DKK21,oster2022approximating,dolgov2022data}.

Our main goal here  is to describe the coupling between our building blocks, explain how they can be  implemented and show our first the numerical results on classical control problems for NS equations. We believe that the simulations presented in the last section illustrate that DP is now feasible from a computational point of view also for fluids and we hope that this can open the way to its application in real industrial applications.\\

The paper is organized as follows.\\
In the second section we will introduce some classical control problems for NS equations and recall the results available in the literature for continuous problems.
Section 3 will be devoted to the presentation of the multilinear approximation and its implementation. In Section 4 we present the TSA and the coupling with the multilinear approximation  In  the last section we present  some numerical experiments  on a number of challenging test problems studied in the literature illustrating the main features of our approach. 

\subsubsection*{Notation and Tensor basics}
In the present work all matrices are represented by large bold-face letters, whereas scalars are given by standard lower-case letters. In the context of model reduction, all matrices with a $\hspace{0.1 cm} \widehat{} \hspace{0.1cm}$ on top, represent low-dimensional quantities.

The Kronecker product of two matrices ${\bm M} \in \RR^{m_1 \times m_2}$ and ${\bm N} \in \RR^{n_1 \times n_2}$ is defined as
$$
{\bm M} \otimes {\bm N} = \begin{pmatrix}
M_{1,1}{\bm N} & \cdots &  M_{1,m_2}{\bm N}\\
\vdots & \ddots & \vdots\\
M_{m_1,1}{\bm N} & \cdots &  M_{m_2,m_2}{\bm N}
\end{pmatrix} \in \RR^{m_1n_1 \times m_2n_2},
$$
and the vec$(\cdot)$ operator maps the entries of a matrix, into a long vector, by stacking the columns of the matrix one after the other. Moreover, we will often make use of the property
\begin{equation*}
\label{kronproperty}
({\bm M} \otimes {\bm N})\mbox{vec}({\bm X}) = \mbox{vec}({\bm N}{\bm X}{\bm M}^{\top}).
\end{equation*}  
Furthermore, the matrix operation $\bm{M} \bullet \bm{N}$, where $\bm{M}$ and $\bm{N}$ have the same size, is known as the Hadamard product, which is an element-element multiplication of the two matrices.

\section{The optimal control problem  for the Navier-Stokes equation and its discretization}\label{sec:ocns}
We introduce some optimal control problems for the Navier-Stokes equation giving also some hints on its numerical solution via finite differences (FD). The approach we present is not limited to this FD approximation and can be extend to other numerical methods as Finite Elements or Finite Volumes since the multilinear discretization described in the next section applies to the semi-discrete system. 

\subsection{The Navier-Stokes dynamical system and its discretization}
Let us assume that the physical domain is a regular bounded connected open set $\Omega$ in $\R^2$ whose boundary will be denoted by $\Gamma$. We denote by $\eta(x)$ the exterior normal vector to a point $x\in\Gamma$. The time variable $t$ will be taken on the interval $(0,T)$ with $T>0$. The standard uncontrolled dynamics will be given by
\begin{equation}
\label{nsgeneral}
\begin{cases}
z_t    -\varepsilon \Delta z + (z\cdot \nabla)z +\nabla p=f  & \hbox{ in } \Omega \times (0,T)\\
\hbox{ div } z =0&  \hbox{ in } \Omega \times (0,T)\\
z=0 & \hbox{ on }\Gamma\times (0,T) \\
z(0)=y_0 &  \hbox{ in }\Omega\\
\end{cases}
\end{equation}

In the sequel we will consider as a model problem the following form for the Navier-Stokes equation
\begin{equation}
\label{nscont}
\begin{cases}
u_t -  \frac{1}{r}\left(u_{xx} + u_{yy}\right)+ uu_x + vu_y  + p_x =f   & \hbox{ in } \Omega \times (0,T)\\
v_t -\frac{1}{r}\left(v_{xx} + v_{yy}\right)+ uv_x +vv_y + p_y = f& \hbox{ in } \Omega \times (0,T) \\
u_x + v_y = 0& \hbox{ in } \Omega \times (0,T)\\
u (t,x)=v(t,x)=0  &\hbox{ on }\Gamma\times (0,T)\\
u(0)=u_0, v(0)=v_0 & \hbox{ in } \Omega \times (0,T) \\
\end{cases}
\end{equation}
where $(x,y)$ are the space coordinates and we set $\varepsilon=1/r$, $z=(u,v)$, where  $u, v: \Omega\times (0,T)\rightarrow \R$, with $\Omega\equiv [0,b_x]\times[0,b_y] \subset \RR^2$ are the velocities to be determined and the pressure $p \in \overline{\Omega} \subset \R^2$ is a Lagrange multiplier introduced to satisfy the incompressibility condition. 
The boundary condition can be chosen according to the problem we want to solve and, as we will see in the next section, can also include some control terms (this will be the case for the boundary control problem). As an example, we can consider the Dirichlet homogeneous boundary condition $z=0$ or the no-slip boundary conditions on each wall, that is
\begin{equation*}
\begin{split}
u(t,x,b_y) &= u_N(x), \quad u(t,x,0) = u_S(x), \quad u(t,0,y) =0, \quad u(t,b_x, y) = 0,\\
v(t,b_x, y) &= v_E(y),\quad v(t,0,y) = v_W(y), \quad v(t,x,0) = 0,  \quad v
(t,x,b_y) =0.\\
\end{split}
\end{equation*}
For the general results on the Navier-Stokes equation we refer to the book \cite{temam2001} whereas for the analysis of some stabilization problems in this framework we refer to \cite{BT2004}.

\subsection{Some control problems for the Navier-Stokes  equation}
The control of non linear dynamical systems over a finite horizon  is usually treated via direct methods based on Pontryagin maximum principle that results in the numerical solution of a two point boundary value problem. This results in an open-loop control and often requires a long work to choose the initial conditions for the state and the co-state since the convergence is local. Moreover, the PMP just gives necessary conditions of optimality and the setting of sufficient conditions is much more technical in this framework.
As we said we are going to present a different approach for some classical control problems that we briefly review here.

\vspace{0.3cm}\noindent
{\em 1. Control acting everywhere in $\Omega$}\\

A first way to introduce a control term is to add a term to the equation, so we write
\begin{equation}\label{eq:controlpb1}
z_t    -\varepsilon \Delta z + (z\cdot \nabla)z +\nabla p+ \sum_{i=1}^m\alpha_i \; \psi_i(x)=f  \qquad \hbox{ in } \Omega \times (0,T).
\end{equation}

In this problem the control is a measurable vector $\alpha:[0,T)\rightarrow A$ where $A$ is a compact subset of $\R^m$ and the functions $\psi_i$ are some predefined shape functions.

\vspace{0.3cm}\noindent
{\em 2. Control acting on a subdomain $\omega\subset\Omega$}\\

Another way to control the NS equation is adding a control term here we apply the control on a subdomain $\omega\subset\Omega$

\begin{equation}\label{eq:controlpb2}
z_t    -\varepsilon \Delta z + (z\cdot \nabla)z + \nabla p=f  + \alpha\; \mathbf{1} _\omega\qquad \hbox{ in } \Omega \times (0,T)
\end{equation}
where $\bold 1_\omega$ is the characteristic function of $\omega$ that typically is a small bounded subdomain of $\Omega$. The control  is a measurable vector $\alpha:[0,T)\rightarrow A$ where $A$ is a compact subset of $\R^m$.

\vspace{0.3cm}\noindent
{\em3.  A boundary control problem on $\gamma\subset \Gamma$}  \\
A third way is  to consider a controlled dynamics  where the control appears in the boundary condition, so we take the dynamics \eqref{nsgeneral} but we modify the boundary Dirichlet condition as 
\begin{eqnarray}\label{eq:controlpb3}
&&z=\alpha  \quad \hbox{ on }\gamma\times (0,T)\\
&&z= g\quad \hbox{ on }(\Gamma\setminus\gamma)\times (0,T)
\end{eqnarray}
and the control must satisfy the compatibility condition 
\begin{equation}
\int_\gamma \alpha  \cdot \eta \;d\gamma=0
\end{equation}
and $\alpha$ will be our control function defined on a small subset $\gamma$ of the boundary $\Gamma$.\\

We will always denote by $y(\cdot, t;\alpha)$ the unique solution of the dynamical system corresponding to the choice of the control $\alpha$. Then we are going to introduce the cost functional to complete the definition of our control problem. 
A general form  of the cost functional we want to minimize is
\begin{equation}
\label{cost}
J(\alpha)= \int_0^T \| y(\cdot, t;\alpha) -\overline y(\cdot, t)\|^2_{L^2(\Omega)} + \gamma \|\alpha\|^2) e^{-\lambda t} dt +\| y(\cdot, T;\alpha) -\overline y(\cdot, T)\|^2_{L^2(\Omega)}
\end{equation}
where the two parameters $\gamma$ and $\lambda$ are positive and the running cost includes the distance of the controlled solution from a reference trajectory $\overline y$  plus a penalization on the control $\alpha$. \\

In the second and third experiment presented in Section
\ref{sec:nens} $\overline y$ will represent the solution of the stationary problem, so the meaning is to stabilize the problem and reach as soon possible the stationary solution.

For the boundary control problem  we will use the same cost functional but $\overline y$ will represent the solution of the problem where we have set a specific control $\overline{\alpha}$.

\subsection{The finite horizon discrete optimal control problem}

Let us recall for reader's convenience the classical Dynamic Programming approach for the {\it finite horizon optimal control problem} for ordinary differential equations, that we use as a model problem. The system is  driven by
\begin{equation}\label{eq}
\left\{ \begin{array}{l}
\dot{y}(s)=f(y(s),\alpha(s),s), \;\; s\in(t,T],\\
y(t)=x\in\R^d,
\end{array} \right.
\end{equation}
and we  denote by $y:[t,T]\rightarrow\R^d$ the solution, by $\alpha:[t,T]\rightarrow\R^m$ the control, by $f:\R^d\times\R^m\times[t,T]\rightarrow\R^d$ the dynamics and by
\[\mathcal{A}=\{\alpha:[t,T]\rightarrow A, \mbox{measurable} \}
\]
the set of admissible controls where $A\subset \R^m$ is a compact set. 
The cost functional for the finite horizon control problem is given by
\begin{equation}\label{cost}
 J_{x,t}(\alpha):=\int_t^T L(y(s,\alpha),\alpha(s),s)e^{-\lambda (s-t)}\, ds+g(y(T))e^{-\lambda (T-t)},
\end{equation}
where $L:\R^d\times\R^m\times [t,T]\rightarrow\R$ is the running cost and $\lambda\geq0$ is the discount factor. The typical assumptions on  the functions $f,L,g$ are:
 \begin{align}
 \begin{aligned}\label{Mf}
|f(x,a,s)| \le M_f,&\quad |L(x,a,s)| \le M_L,\quad |g(x)| \le M_g,  \quad \forall\, x \in \mathbb{R}^d, a \in A \subset \mathbb{R}^m, s \in [t,T], 
\end{aligned}
\end{align}
the functions $f, L$ and $g$ are Lipschitz-continuous with respect to the first variable
\begin{align}
\begin{aligned}\label{Lf}
&|w(x,a,s)-w(y,a,s)| \le L_w |x-y|,\hbox{  for } w=f,L \hbox{  and }\;\forall \, x,y \in \mathbb{R}^d, a \in A \subset \mathbb{R}^m, s \in [t,T], 
\end{aligned}
\end{align}
$$
|g(x)-g(y)| \le L_g |x-y|, \qquad\forall \, x,y \in \mathbb{R}^d.
$$
%

Note that these assumptions guarantee uniqueness for the trajectory $y(t)$ by the Carath\'eodory theorem  (we refer to e.g. \cite{BCD97} for a precise statement).

The goal is to find a state-feedback control law $\alpha(t)=\Phi(y(t),t),$ in terms of the state equation $y(t),$ where $\Phi$ is the feedback map. To derive optimality conditions we use the well-known Dynamic Programming Principle (DPP) due to Bellman. We first define the value function for an initial condition $(x,t)\in\R^d\times [t,T]$:
\begin{equation}
v(x,t):=\inf\limits_{\alpha\in\mathcal{A}} J_{x,t}(\alpha)
\label{value_fun}
\end{equation}
A classical result (see \cite{BCD97} shows that under our assumptions the value function for the finite horizon problem is the unique viscosity solution of the following Hamilton-Jacobi-Bellman equation  

\begin{equation}\label{HJB}
\left\{
\begin{array}{ll} 
-\dfrac{\partial v}{\partial s}(x,s) +\lambda v(x,s)+ \max\limits_{a\in A }\left\{-L(x, a,s)- \nabla v(x,s) \cdot f(x,a,s)\right\} = 0, & \hbox { for } x\in\R^d \, s\in [t,T) \\
v(x,T) = g(x) & \hbox{ for } x\in \R^d
\end{array}
\right.
\end{equation}
Once the value function is known, by e.g. \eqref{HJB}, then it is possible to  compute the optimal feedback control as:
\begin{equation}\label{feedback}
\alpha^*(t):=  \argmax_{a \in A }\left\{-L(x,a,t)- \nabla v(x,t) \cdot f(x,a,t)\right\}. 
\end{equation}
and this is one of the most important features of the DP approach to control problems.
From the numerical approximation of the feedback we can apply \eqref{feedback} replacing 
the continue value function $v$ by our numerical approximation and its gradient by discrete gradients (see \cite{FF14} for more details on this point).

\section{Multilinear approximation of the Navier-Stokes equation}

In this section we present the matrix-oriented discretization of the NS equation used to set-up the finite dimensional optimal control problem. The technique presented here allows us to significantly speed up the 
numerical time integration; see the first test presented in Section \ref{sec:nens}. More details on the procedure presented here can be found in \cite{kirsten_phd}.

\subsection{Matrix-Oriented discretization and 2S-POD-DEIM for general PDEs}
\label{2spod}

Consider a semilinear evolutive PDE of the form
\begin{equation} \label{nonlinpdeintro}
u_{t} = {\cal L}(u) + f(u,t) , \quad u = u({\bm x}, t) \quad \mbox{with} \,\,
 {\bm x} \in \Omega \subset \mathbb{R}^2,\,\,\, t\in \mathcal{T},
\end{equation} 
with suitable boundary conditions. We assume that the differential operator ${\cal L}$ is linear in $u$ with separable coefficients, typically a second order
operator in the space variables, while 
{$f: S \times \mathcal{T} \rightarrow \RR$} is a nonlinear function,
 where $S$ is an appropriate space with $u\in S$, and $\mathcal{T}$ is the timespan. Under these assumptions, if ${\cal L}$ is discretized by means of a tensor basis, such as finite differences on rectangular domain domains, certain finite element methods and certain spectral methods, then the physical domain can be mapped to a reference hypercubic domain. More precisely, if we consider finite differences on a rectangular domain, then 
 $$
 {\cal L} = \bm{I}_{n_y} \otimes \bm{A}_1 + \bm{A}_2^{\top} \otimes I_{n_x},
 $$
 where $\bm{A}_1 \in \mathbb{R}^{n_x \times n_x}$ and $\bm{A}_2 \in \mathbb{R}^{n_y \times n_y}$ are matrices containing the coefficients for the derivatives and $n_x$ and $n_y$ are the number of discretization nodes in the $x-$ and $y-$ directions respectively. As a result if we define ${\bm U}(t) \in \mathbb{R}^{n_x \times n_y}$ as a matrix containing an approximation to the solution $u(t)$ at each discretization node, then the discrete version of \eqref{nonlinpdeintro} can be expressed in matrix form as
\begin{equation}
 \bm{\dot{U}}(t) = \bm{A}_1\bm{U}(t) + \bm{U}(t)A_2 + \bm{F}(\bm{U}, t).
 \label{matdiff}   
\end{equation}
In addition to a better structural interpretation of the discrete quantities, this formulation can also lead to reduced memory requirements and computational costs. A summary of this and related matrix-oriented procedures can be found in \cite{Simoncini2017, palitta2016}. Furthermore, standard numerical integration schemes, such as semi-implicit schemes and exponential integrators, can be performed directly in matrix form to approximate the solution of \eqref{matdiff} throughout the timespan \cite{dautilia20, kirsten_phd}.

As it is well-known in the vector formulation of discretized PDEs, the discrete matrices are often very large and sparse and require a large computational effort to solve the resulting linear systems at each time step. To this end model order reduction techniques such as POD \cite{hinze2005proper, BCOW.17} and DEIM \cite{chaturantabut2010nonlinear} have been successfully applied to reduce the complexity of solving several linear systems throughout the timespan.In \cite{Kirsten.Simoncini.arxiv2020, kirsten.21} the POD and DEIM methods have been extended so that they can be applied directly to the matrix differential equation \eqref{matdiff}, without requiring any mapping from matrices to vectors.

In short,consider a set of $n_s$ time-dependent snapshot solutions and nonlinear snapshots of \eqref{matdiff}, given by $\{\bm{U}^i\}_{i = 1}^{n_s}$ and  $\{{\bm F}(\bm{U}^i)\}_{i = 1}^{n_s}$ respectively. The 2S-POD algorithm from \cite{Kirsten.Simoncini.arxiv2020} is applied to the set of snapshot solutions to form two tall matrices $\bm{U}_{\ell} \in \mathbb{R}^{n_x \times k_{\ell}}$ and $\bm{U}_{r} \in \mathbb{R}^{n_y \times k_r}$ ($k_\ell$, $k_r \ll n$)  with orthornormal columns. The parameters $k_\ell$ and $k_r$ refer to the number of selected dominant singular values such that the projection error is bounded by a prescribed tolerance $tol$. The span of these columns respectively approximate the row and column space of the snapshot solutions. To this end, we approximate the solution $\bm{U}(t)$ of \eqref{matdiff} by $\bm{U}(t) \approx \bm{U}_{\ell}\widehat{\bm{U}}(t)\bm{U}_{r}^{\top} =: \widetilde{\bm{U}}(t)$, for $t \in \mathcal{T}$,  where $\widehat{\bm{U}}(t) \in \mathbb{R}^{k_\ell \times k_r}$ satisfies the following low-dimensional matrix ODE:
\begin{equation}
\begin{cases}
\dot{\widehat{\bm{U}}}(t) = \widehat{\bm{A}}_1\widehat{\bm{U}}(t) + \widehat{\bm{U}}(t) \widehat{\bm{A}}_2 + \widehat{F}(\widetilde{\bm{U}}(t), t), & t \in \mathcal{T},
\\
\widehat{\bm{U}}(0)=\widehat{\bm{U}}_0 = \bm{U}_{\ell}^{\top}\bm{U}_0\bm{U}_{r},
\end{cases}
\label{smallmatrixOde}
\end{equation}
where $\widehat{\bm{A}}_1 = \bm{U}_{\ell}^{\top} \bm{A}_1\bm{U}_{r}$ and $\widehat{\bm{A}}_2 = \bm{U}_{\ell}^{\top} \bm{A}_2\bm{U}_{r}$ and $\widehat{F}(\widetilde{\bm{U}}(t), t) = \bm{U}_{\ell}^{\top}{\bm F}(\bm{U}_{\ell}\widehat{\bm{U}}(t)\bm{U}_{r}^{\top}, t)\bm{U}_{r}$. Despite the fact that $\widehat{F}(\widetilde{\bm{U}}(t), t)$ is considered a low-dimensional quantity, the calculation of it still results in a computational bottleneck, since the nonlinear function first needs to be evaluated at all the entries of $\bm{U}_{\ell}\widehat{\bm{U}}(s)\bm{U}_{r}^{\top}\in \mathbb{R}^{n_x \times n_y}$ before it is projected onto the low-dimensional space. To overcome this bottleneck we apply the 2S-DEIM method from \cite{Kirsten.Simoncini.arxiv2020}.

We consider the set of nonlinear snapshots $\{{\bm F}({\bm U}^i)\}_{i = 1}^{n_s}$ and use the 2S-POD algorithm to form two tall matrices $\bm{\Phi}_{\ell} \in \mathbb{R}^{n_x \times p_1}$ and $\bm{\Phi}_{r}  \in \mathbb{R}^{n_y \times p_2}$ ($p_i \ll n$)  with orthornormal columns. We aim to approximate the nonlinear term by far smaller matrices, that is ${\bm F}(\bm{U}_{\ell}\widehat{\bm{U}}(s)\bm{U}_{r}^{\top}) \approx \bm{\Phi}_{\ell}\widehat{\bm F}(t)\bm{\Phi}_{r}^{\top}$, where $\widehat{\bm F}(t) \in \mathbb{R}^{p_1 \times p_2}$ is a matrix of time-dependent coefficients. This leads to a 2S-DEIM approximation of the form:
\begin{equation}
\widehat{F}(\widetilde{\bm U}(t), t) \approx \bm{\Phi}_{\ell}(\bm{D}_{\ell}^\top \bm{\Phi}_{\ell})^{-1}\bm{D}_{\ell}^\top{\bm F}(\widetilde{\bm U}(t), t)\bm{D}_{r}(\bm{\Phi}_{r}^\top \bm{D}_{r})^{-1}\bm{\Phi}_{r}^\top,
\label{DEIM}
\end{equation}
where $\bm{D}_{\ell} \in \mathbb{R}^{n_x \times p_1}$ and $\bm{D}_{r} \in \mathbb{R}^{n_y \times p_2}$ respectively contain a subset of $p_1$ and $p_2$ columns of the $n_x \times n_y$ identity matrix. The indices at which these columns are selected is determined by respectively applying the Q-DEIM algorithm from \cite{drmac2016} to the matrices $\bm{\Phi}_{\ell}$ and $\bm{\Phi}_{r}$. In the elegant case where the nonlinear function is evaluated elementwise at the indices of $\widetilde{\bm U}(t)$, the respective interpolation indices can be selected by taking $\bm{D}_{\ell}$ and $\bm{D}_{r}$ inside the nonlinear function such that

\begin{equation}
\widehat{F}(\widetilde{\bm U}(t), t) \approx \bm{U}_\ell^\top \bm{\Phi}_{\ell}(\bm{D}_{\ell}^\top \bm{\Phi}_{\ell})^{-1}{\bm F}(\bm{D}_{\ell}^\top \widetilde{\bm U}(t)\bm{D}_{r},t)(\bm{\Phi}_{r}^\top \bm{D}_{r})^{-1}\bm{\Phi}_{r}^\top \bm{U}_r,
\label{fapprox}
\end{equation}
can be completely evaluated in low-dimension. In the case that the nonlinear term is not evaluated elementwise, more complex techniques may be required \cite{chaturantabut2010nonlinear}. This situation is also encountered with the NS equation and will be discussed in the following section.

In what follows we aim to extend these matrix-oriented discretization, integration and model reduction strategies to the setting of the NS equation.

\subsection{The NS equation in full dimension}
\label{SecNSfull}
The method considered for the time and space discretization of \eqref{nscont} is finite differences on a staggered grid. A discussion of the scheme and a Matlab implementation in the vector setting can be respectively found  in \cite{strang07} and \cite{seibold08}. Here we aim to take explicit advantage of the rectangular domain, to directly treat the equation in matrix form, both for the reduction and integration phases of the method.

For the space discretization, we consider $n_x$ gridpoints in the $x-$direction and $n_y$ gridpoints in the $y-$direction. For the staggered grid, the velocities $u$ are placed on the vertical cell interfaces, $v$ on the horizontal cell interfaces and the pressure $p$ in the centre of the cells. That is, the discretized quantities are given by the matrices ${\bm U}(t) \in \RR^{n_x-1 \times n_y}$, ${\bm V}(t) \in \RR^{n_x \times n_y-1}$ and ${\bm P}(t) \in \RR^{n_x \times n_y}$. Given $\ast \in \{U,V\}$, we consider the matrices ${\bm A}_{1, \ast}$ and ${\bm A}_{2, \ast}$ with corresponding dimensions to respectively contain the coefficients for the second derivative in the $x-$ and $y-$ directions and the matrices ${\bm B}_{1, \ast}$ and ${\bm B}_{2, \ast}$ that of the first derivatives. The discrete version of \eqref{nscont} is then given by
\begin{equation}
\label{nsdisc}
\begin{cases}
\dot{\bm U}-{\bm A}_{1, U}{\bm U} + {\bm U}{\bm A}_{2, U}^{\top} + {\bm B}_{1,U}^\top {\bm P} + {\bm F}_U({\bm U}, {\bm V}, t)  &= 0 \\
\dot{\bm V}-{\bm A}_{1, V}{\bm V} + {\bm V}{\bm A}_{2, V}^{\top}  +  {\bm P}{\bm B}_{2,V} + {\bm F}_V({\bm U}, {\bm V}, t) &= 0 \\
{\bm B}_{1, U}{\bm U} + {\bm V}{\bm B}_{2,V}^{\top} &= 0,
\end{cases}
\end{equation}
where ${\bm F}_U({\bm U}, {\bm V}, t) = {\bm B}_{1, U}{\bm U}\bullet{\bm U} + {\bm U}{\bm B}_{2,U}^{\top}\bullet{\bm V}$ and ${\bm F}_V({\bm U}, {\bm V}, t) = {\bm B}_{1, V}{\bm V}\bullet{\bm U} + {\bm V}{\bm B}_{2,V}^{\top}\bullet{\bm V}$.
For the time discretization we consider a simple semi-implicit Euler scheme, so that the viscosity terms are treated implicitly, the nonlinear terms explicitly, and the pressure term is treated implicitly via a Chorin Projection scheme (see \cite{chorin68}). That is, at each time iteration the approximations ${\bm U}^{(j+1)} \approx {\bm U}(t_{j+1})$
and ${\bm V}^{(j+1)} \approx {\bm V}(t_{j+1})$ are determined by solving the Sylvester equations
{\small \begin{equation}
\label{nstimedisc}
\begin{cases}
(\bm I- \Delta t{\bm A}_{1, U}){\bm U}^{(j+1)} + {\bm U}^{(j+1)}(-\Delta t{\bm A}_{2, U}^{\top}) = {\bm U}^{(j)} - {\bm B}_{1,U}^\top {\bm P}^{(j+1)} - \Delta t{\bm F}_{U}({\bm U}^{(j)}, {\bm V}^{(j)}, t) \\
(\bm I- \Delta t{\bm A}_{1, V}){\bm V}^{(j+1)} + {\bm V}^{(j+1)}(-\Delta t{\bm A}_{2, V}^{\top}) = {\bm V}^{(j)} -  {\bm P}^{(j+1)}{\bm B}_{2,V} - \Delta t{\bm F}_{V}({\bm U}^{(j)}, {\bm V}^{(j)}, t).
\end{cases}
\end{equation}}
Keeping the pressure term at the next time step on the right hand side of the Sylvester equations is a slight abuse of notation. In fact, the pressure is determined by a pressure correction to enforce the incompressibility.

More precisely, consider the implicit time discretization of the pressure, that is ${\bm U}^{(j+1)} - {\bm U}^{(j)} = \Delta t {\bm B}_{1,U}^\top {\bm P}^{(j+1)}$ and ${\bm V}^{(j+1)} - {\bm V}^{(j)} = \Delta t {\bm P}^{(j+1)} {\bm B}_{2,V}$. If we multiply the first equation from the left by ${\bm B}_{1,U}$ and the second from the right by ${\bm B}_{2,V}^\top$ adding the two equations together, we obtain a Sylvester equation of the form 
\begin{equation}
{\bm A}_{1,U}{\bm P}^{(j+1)} + {\bm P}^{(j+1)}{\bm A}_{2,V}^{\top} = {\bm B}_{1,U}{\bm U}^{(j)} + {\bm V}^{(j)}{\bm B}_{2,V}^\top,
\label{pres_corr}
\end{equation}
to be solved for ${\bm P}^{(j+1)}$. This equation is obtained by enforcing the incompressibility such that ${\bm B}_{1,U}{\bm U}^{(j+1)} + {\bm V}^{(j+1)}{\bm B}_{2,V}^\top = 0$. Determining the pressure correction does not cause any problems with respect to the staggered grid, since the divergence of the velocity lies in the cell centres, similar to the pressure. Determining the nonlinear terms is, however, more complicated.

More precisely, due to the staggered grid, the nodes of the of $\bm{U}$ and $\bm{V}$ are located at different positions so that, for example, the product ${\bm U} \bullet {\bm V}$ is not defined. This is circumvented by means of interpolation, for which we refer to [\cite{seibold08}, section 5] for details. To this end, following [\cite{seibold08}, section 5] to incorporate the boundary conditions in the nonlinear terms, we define $\overline{\bm U} \in \mathbb{R}^{(n_x+1) \times n_y}$ as the matrix $\bm{U}$ padded with boundary conditions in the top and bottom rows. Similarly, $\overline{\bm V} \in \mathbb{R}^{n_x \times n_y + 1}$ is padded with boundary conditions in the first and last columns. Then, defining $\bm{C} \in \mathbb{R}^{(n_x - 1) \times (n_x + 1)}$ as a matrix with $1/2$ on the main and upper diagonal as an averaging matrix, the nonlinear term, evaluated on the staggered grid, can be expressed in fully matricial form as 
\begin{equation}
{\bm F}_{U}({\bm U}, {\bm V}, t) = {\bm B}_{1,U}^{\top} \left((\bm{C}^{\top} \overline{\bm U})^2-\gamma |\bm{C}^{\top} \overline{\bm U}|\bullet \left(\frac{h_x}{2} {\bm B}_{1,U}^{\top}\overline{\bm U} \right) \right) +\overline{\bm B}_{1,U}^{\top} \left((\overline{\bm U} \overline{\bm{C}} ) \bullet (\bm{C}^{\top} \overline{\bm V})-\gamma |(\overline{\bm U} \overline{\bm{C}} |\bullet \left(\frac{h_x}{2} \overline{{\bm B}}_{2,V}^{\top}\overline{\bm V} \right) \right),
\label{nonlinearfunction}
\end{equation}
where $\frac{h_x}{2} {\bm B}_{1,U}^{\top}$ is responsible for differencing and cofficient matrices containing a superscripted bar are merely conforming to the dimension of $\bm{\overline{U}}$ or $\bm{\overline{V}}$. A similar form can be derived for ${\bm F}_{V}({\bm U}, {\bm V}, t)$.

A summary of the procedure for solving \eqref{nscont} on a staggered grid in fully matricial form can be found in Algorithm \ref{alg_full}. The most computationally expensive step is the solution of three Sylvester equations (one at Step 3, and two at step 5) at each timestep. However, the coefficient matrices remain constant throughout the timespan. To this end, an a-priori eigenvalue decomposition of the six coefficient matrices can be performed so that the Sylvester equations can be solved by using only substitution and matrix-matrix multiplication. We refer the reader to \cite{Simoncini2017, dautilia20, kirsten_phd} for further details.

\begin{algorithm}[H]
\caption{Full model Algorithm}
\begin{algorithmic}[1]
\State Choose the initial conditions $({\bm U}(0)$, ${\bm V}(0))$ and the number of time steps $n_t$ 
\For{$i=0,\ldots, n_t-1$}
\State Solve \eqref{pres_corr} for ${\bm P}^{(i+1)}$ 
\State Compute  ${\bm F}_{U}({\bm U}, {\bm V}, t)$ and  ${\bm F}_{V}({\bm U}, {\bm V}, t)$
\State Solve \eqref{nsdisc}
\EndFor
\end{algorithmic}
\label{alg_full}
\end{algorithm}

\subsection{The 2S-POD-DEIM reduced NS equation}
For the reduced model we consider the 2S-POD-DEIM model reduction procedure for systems of matrix differential equations from \cite{kirsten2019,kirsten.21}, discussed above, to reduce \eqref{nsdisc} in a fully matricial way. To this end we consider $n_s$ snapshots of the full dimensional solutions $\bm{U}_i, \bm{V}_i$ and $\bm{P}_i$, $i = 1,2,\ldots,n_s$, in order to construct the low dimensional, orthonormal basis matrices $\bm{U}_\ast \in \RR^{n_x - 1 \times k_{1,\ast}}$, $\bm{V}_\ast \in \RR^{n_y - 1 \times k_{2,\ast}}$ and $\bm{P}_\ast \in \RR^{n_x \times k_{3,\ast}}$, where $\ast = \{\ell,r\}$. This leads to the approximations
$$
    {\bm U} \approx {\bm U}_\ell \widehat{\bm U} {\bm U}_r^{\top} =: \widetilde{\bm U}, \qquad  {\bm V} \approx {\bm V}_\ell \widehat{\bm V} {\bm V}_r^{\top} =: \widetilde{\bm V} \quad \mbox{and} \quad {\bm P} \approx {\bm P}_\ell
    \widehat{\bm P} {\bm P}_r^{\top} =: \widetilde{\bm P}.
$$
\subsubsection{Solving for $\widehat{\bm U}^{(j+1)}$ and $\widehat{\bm V}^{(j+1)}$}
Substituting the above approximations into \eqref{nsdisc} yields a reduced Navier-Stokes equations, where the reduced solutions $\widehat{\bm U}^{j+1}$ and $\widehat{\bm V}^{j+1}$ are determined by solving the (reduced) coupled Sylvester equations
{\small \begin{equation}
\label{nstimedisc}
\begin{cases}
(\bm I- \Delta t\widehat{\bm A}_{1, U})\widehat{\bm U}^{(j+1)} + \widehat{\bm U}^{(j+1)}(-\Delta t\widehat{\bm A}_{2, U}^{\top}) = \widehat{\bm U}^{(j)} - \bm{U}_\ell^\top{\bm B}_{1,U}^\top \bm{P}_\ell\widehat{\bm P}^{(j+1)}\bm{P}_r^\top\bm{U}_r - \Delta t\widehat{\bm F}_{U}(\widetilde{\bm U}^{(j)}, \widetilde{\bm V}^{(j)}, t) \\
(\bm I- \Delta t\widehat{\bm A}_{1, V})\widehat{\bm V}^{(j+1)} + \widehat{\bm V}^{(j+1)}(-\Delta t\widehat{\bm A}_{2, V}^{\top}) = \widehat{\bm V}^{(j)} -  \bm{V}_\ell^\top\bm{P}_\ell\widehat{\bm P}^{(j+1)}\bm{P}_r^\top{\bm B}_{2,V}\bm{V}_r - \Delta t\widehat{\bm F}_{V}(\widetilde{\bm U}^{(j)}, \widetilde{\bm V}^{(j)}, t).
\end{cases}
\end{equation}}
Here we have used the same numerical integration scheme as for the full-dimensional equation, and all the matrices that have a$\hspace{0.1cm} \widehat{} \hspace{0.1cm}$ on the top are left and right projections of the original coefficient matrices onto the relevant subspaces, and hence they are all low-dimensional. Furthermore, notice that the terms multiplying the pressure term from the left and right in both equations are low-dimensional and time-independent, hence they can be stored offline. Consequently, the remaining challenges lie in determining the pressure term $\widehat{\bm P}^{(j+1)}$ and evaluating the nonlinear functions in low-dimension at each timestep. 

\subsubsection{Solving for $\widehat{\bm P}^{(j+1)}$}

 The pressure correction step requires the solution of the Sylvester equation \eqref{pres_corr}. Inserting the approximation $\widetilde{\bm{P}}^{j+1}$ into \eqref{pres_corr} yields the low-dimensional Sylvester equation 
\begin{equation}
\bm{P}_\ell^\top{\bm A}_{1,U}\bm{P}_\ell\widehat{\bm P}^{(j+1)} + \widehat{\bm P}^{(j+1)}\bm{P}_r^{\top}{\bm A}_{2,V}^{\top}\bm{P}_r = {\bm P}_\ell^{\top} {\bm B}_{1,U}\overline{\bm{U}_\ell\widehat{\bm U}^{(j)}\bm{U}_r^\top}{\bm P}_r  + {\bm P}_\ell^{\top} \overline{\bm{V}_\ell\widehat{\bm V}^{(j)}\bm{V}_r^\top}{\bm B}_{2,V}^\top {\bm P}_r.
\label{pres_corr_red}
\end{equation}
 Here, $\overline{\bm{U}_\ell\widehat{\bm U}^{(j)}\bm{U}_r^\top}$ and $\overline{\bm{V}_\ell\widehat{\bm V}^{(j)}\bm{V}_r^\top}$ represent the lifted quantities padded with boundary conditions, as discussed before. As a result, the left-hand side of this Sylvester equation consists of only small matrices, but the right hand side, on the other hand, requires some more attention to avoid recomputing large matrices, due to the fact that the lifted quantities are padded by boundary conditions. To this end, by taking advantage of the fact that ${\bm B}_{1,U}$ and ${\bm B}_{2,V}$ only account for the differentiation, the first and last rows of $\bm{P}_\ell$ and $\bm{P}_r$ can be manipulated in such a way that they only act on the boundary conditions so that the internal blocks of $\bm{P}_\ell$ and ${\bm B}_{1,U}$ ($\bm{P}_r$ and ${\bm B}_{2,V}$) can be multiplied with $\bm{U}_\ell$ ($\bm{U}_r$) in order to form low-dimensional matrices that can be stored offline. A similar manipulation is done for the products between $\bm{U}_r^\top$ and $\bm{P}_r$ ($\bm{P}_\ell^\top$ and $\bm{V}_\ell$) such that only low-dimensional operations need to occur online.

\subsubsection{Evaluating $\widehat{F}_U(\widetilde{\bm U},\widetilde{\bm V},t)$ and $\widehat{F}_V(\widetilde{\bm U},\widetilde{\bm V},t)$ with DEIM}

Following the procedure in \cite{Kirsten.Simoncini.arxiv2020}, we consider $n_S$ snapshots of the nonlinear functions $\bm{F}_{U}$ and $\bm{F}_V$ to construct the low-dimensional, orthonormal matrices $\bm{\Phi}_{\bm{U}, \ast} \in \mathbb{R}^{n_x -1 \times p_{1,\ast}}$ and $\bm{\Phi}_{\bm{V}, \ast} \in \mathbb{R}^{n_y -1 \times p_{2,\ast}}$, $\ast = \{\ell,r\}$ used for the reduction of the nonlinear function by 2S-DEIM. If we define the orthonormal matrices $\bm{D}_{\bm{U}, \ast} \in \mathbb{R}^{n_x -1 \times p_{1,\ast}}$  and $\bm{D}_{\bm{V}, \ast} \in \mathbb{R}^{n_y -1 \times p_{2,\ast}}$ as matrices with a subset of columns of the identity matrix, then the 2S-DEIM approximation of the nonlinear terms is given by  
\begin{equation}
\widehat{F}_{U}(\widetilde{\bm U},\widetilde{\bm V},t) \approx \bm{U}_\ell^{\top} \bm{\Phi}_{\bm{U}, \ell} \left(\bm{D}_{\bm{U}, \ell}^{\top} \bm{\Phi}_{\bm{U}, \ell} \right)^{-1} \bm{D}_{\bm{U}, \ell}^{\top} F_U(\widetilde{\bm U},\widetilde{\bm V},t) \bm{D}_{\bm{U}, r} \left( \bm{\Phi}_{\bm{U}, r}^{\top} \bm{D}_{\bm{U}, r} \right)^{-1}  \bm{\Phi}_{\bm{U}, r}^{\top} \bm{U}_r,
\label{term_deim}
\end{equation}
and similar for $\widehat{F}_{V}(\widetilde{\bm U},\widetilde{\bm V},t)$. The matrices $\bm{U}_\ell^{\top} \bm{\Phi}_{\bm{U}, \ell} \left(\bm{D}_{\bm{U}, \ell}^{\top} \bm{\Phi}_{\bm{U}, \ell} \right)^{-1} $ and $\bm{U}_r^{\top} \bm{\Phi}_{\bm{U}, r} \left(\bm{D}_{\bm{U}, r}^{\top} \bm{\Phi}_{\bm{U}, r} \right)^{-1} $ are low-dimensional and can be stored offline, however in this setting it is particularly challenging to evaluate the term $\bm{D}_{\bm{U}, \ell}^{\top} F_U(\widetilde{\bm U},\widetilde{\bm V},t) \bm{D}_{\bm{U}, r}$ without first lifting and evaluating $F_U(\widetilde{\bm U},\widetilde{\bm V},t)$ in full dimension. In what follows we briefly discuss how this is achieved. The same idea follows for $F_V(\widetilde{\bm U},\widetilde{\bm V},t)$.

We want to determine \eqref{term_deim} by using only small matrices. From \eqref{nonlinearfunction}, it can be seen that $F_U$ consists of two terms summed together. Therefore:
$$
\bm{D}_{\bm{U}, \ell}^{\top} F_U(\widetilde{\bm U},\widetilde{\bm V},t) \bm{D}_{\bm{U}, r} = \bm{D}_{\bm{U}, \ell}^{\top}F_{U,1}(\widetilde{\bm U},\widetilde{\bm V},t) \bm{D}_{\bm{U}, r}+ \bm{D}_{\bm{U}, \ell}^{\top} F_{U,2}(\widetilde{\bm U},\widetilde{\bm V},t) \bm{D}_{\bm{U}, r}.
$$

In the following result we illustrate how the first of the two terms can be evaluated in low dimension. A similar strategy is used for the second term, but for the sake of presentation we omit the details.

\begin{prop}
The term $\bm{D}_{\bm{U}, \ell}^{\top}F_{U,1}(\widetilde{\bm U},\widetilde{\bm V},t) \bm{D}_{\bm{U}, r}$ can be evaluated completely in low-dimension, independent of the full dimensions $n_x$ and $n_y$, at each time step.
\end{prop}

{\em Proof.}

From the definition of the full-dimensional nonlinear term, it can be seen that 
$$
\bm{D}_{\bm{U}, \ell}^{\top}F_{U,1}(\widetilde{\bm U},\widetilde{\bm V},t) \bm{D}_{\bm{U}, r} = \bm{D}_{\bm{U}, \ell}^{\top} {\bm B}_{1,U}^{\top} \left((\bm{C}^{\top} \widetilde{\bm U})^2-\gamma |\bm{C}^{\top} \widetilde{\bm U}|\bullet \left(\frac{h_x}{2} {\bm B}_{1,U}^{\top}\widetilde{\bm U} \right) \right)\bm{D}_{\bm{U}, r}.
$$
As a result, 
$$
\bm{D}_{\bm{U}, \ell}^{\top} {\bm B}_{1,U}^{\top} \left((\bm{C}^{\top} \widetilde{\bm U})^2-\gamma |\bm{C}^{\top} \widetilde{\bm U}|\bullet \left(\frac{h_x}{2} {\bm B}_{1,U}^{\top}\widetilde{\bm U} \right) \right)\bm{D}_{\bm{U}, r} 
$$
$$
= \bm{D}_{\bm{U}, \ell}^{\top} {\bm B}_{1,U}^{\top} \left(\bm{C}^{\top} \widetilde{\bm U}\bullet \bm{C}^{\top} \widetilde{\bm U}\right) \bm{D}_{\bm{U}, r} - \bm{D}_{\bm{U}, \ell}^{\top} {\bm B}_{1,U}^{\top} \gamma |\bm{C}^{\top} \widetilde{\bm U}|\bullet \left(\frac{h_x}{2} {\bm B}_{1,U}^{\top}\widetilde{\bm U} \right) \bm{D}_{\bm{U}, r}.
$$
Once again we look at the two terms on the right-hand side separately. More precisely, considering the first term, the role of $\bm{D}_{\bm{U}, r}$ is to select columns after the scalar product. Therefore, it can be taken inside of the scalar product, such that
$$
\bm{D}_{\bm{U}, \ell}^{\top} {\bm B}_{1,U}^{\top} \left(\bm{C}^{\top} \widetilde{\bm U}\bullet \bm{C}^{\top} \widetilde{\bm U}\right) \bm{D}_{\bm{U}, r} = \bm{D}_{\bm{U}, \ell}^{\top} {\bm B}_{1,U}^{\top} \left(\bm{C}^{\top} {\bm U}_\ell  \widehat{\bm U}(t) {\bm U}_r^\top \bullet \bm{C}^{\top} {\bm U}_\ell  \widehat{\bm U}(t) {\bm U}_r^\top    \right) \bm{D}_{\bm{U}, r}
$$
$$
= \bm{D}_{\bm{U}, \ell}^{\top} {\bm B}_{1,U}^{\top} \left(\bm{C}^{\top} {\bm U}_\ell  \widehat{\bm U}(t) {\bm U}_r^\top \bm{D}_{\bm{U}, r} \bullet \bm{C}^{\top} {\bm U}_\ell  \widehat{\bm U}(t) {\bm U}_r^\top \bm{D}_{\bm{U}, r}    \right)
$$
The term ${\bm U}_r \bm{D}_{\bm{U}, r}$ is small and can be saved offline.
The role of the term $\bm{D}_{\bm{U}, \ell}^{\top} {\bm B}_{1,U}^{\top}$ is to select at which rows the derivatives is taken after the scalar product. Therefore for each row $e_i$ selected we need to compute $\frac{e_{i+1}-e_i}{h_x}$. This means we only need the scalar product at rows $i$ and $i+1$.
Hence,
$$
\bm{D}_{\bm{U}, \ell}^{\top} {\bm B}_{1,U}^{\top} \left(\bm{C}^{\top} {\bm U}_\ell  \widehat{\bm U}(t) {\bm U}_r^\top \bm{D}_{\bm{U}, r} \bullet \bm{C}^{\top} {\bm U}_\ell  \widehat{\bm U}(t) {\bm U}_r^\top \bm{D}_{\bm{U}, r}    \right)
$$
$$
= \frac{1}{h_x} \left( (\bm{D}_{\bm{U}, \ell}^+)^{\top} \bm{C}^{\top} {\bm U}_\ell  \widehat{\bm U}(t) {\bm U}_r^\top \bm{D}_{\bm{U}, r} \bullet (\bm{D}_{\bm{U}, \ell}^+)^{\top} \bm{C}^{\top} {\bm U}_\ell  \widehat{\bm U}(t) {\bm U}_r^\top \bm{D}_{\bm{U}, r}    \right)
$$
$$
- \frac{1}{h_x} \left( \bm{D}_{\bm{U}, \ell}^{\top} \bm{C}^{\top} {\bm U}_\ell  \widehat{\bm U}(t) {\bm U}_r^\top \bm{D}_{\bm{U}, r} \bullet \bm{D}_{\bm{U}, \ell}^{\top} \bm{C}^{\top} {\bm U}_\ell  \widehat{\bm U}(t) {\bm U}_r^\top \bm{D}_{\bm{U}, r}    \right) 
$$
can be computed via only low-dimensional evaluations, since the only time-dependent term is $\widehat{\bm{U}}(t)$ and all other matrix products result in low-dimensional coefficient matrices that can be stored offline. The matrix $\bm{D}_{\bm{U}, \ell}^+$ contains the Deim indices shifted by $+1$. The same idea works for the second term  
$$
\bm{D}_{\bm{U}, \ell}^{\top} {\bm B}_{1,U}^{\top} \gamma |\bm{C}^{\top} \widetilde{\bm U}|\bullet \left(\frac{h_x}{2} {\bm B}_{1,U}^{\top}\widetilde{\bm U} \right) \bm{D}_{\bm{U}, r}
$$
which can be expressed as
{\small
$$
\frac{1}{h_x} \gamma \left| (\bm{D}_{\bm{U}, \ell}^+)^{\top} \bm{C}^{\top} {\bm U}_\ell  \widehat{\bm U}(t) {\bm U}_r^\top \bm{D}_{\bm{U}, r} \right| \bullet \left( \frac{h_x}{2} (\bm{D}_{\bm{U}, \ell}^+)^{\top} {\bm B}_{1,U}^{\top} {\bm U}_\ell  \widehat{\bm U}(t) {\bm U}_r^\top \bm{D}_{\bm{U}, r}   \right)
$$
$$
- \frac{1}{h_x} \gamma \left| \bm{D}_{\bm{U}, \ell}^{\top} \bm{C}^{\top} {\bm U}_\ell  \widehat{\bm U}(t) {\bm U}_r^\top \bm{D}_{\bm{U}, r} \right| \bullet \left( \frac{h_x}{2} \bm{D}_{\bm{U}, \ell}^{\top} {\bm B}_{1,U}^{\top} {\bm U}_\ell  \widehat{\bm U}(t) {\bm U}_r^\top \bm{D}_{\bm{U}, r}   \right).
$$}
Once again all coefficient matrices are low-dimensional and stored offline, so that only low-dimensional matrix multiplications are necessary to obtain the term $\bm{D}_{\bm{U}, \ell}^{\top} F_{U,1}(\widetilde{\bm U},\widetilde{\bm V},t) \bm{D}_{\bm{U}, r}$. This completes the proof.

A brief summary of the reduced model phase has been sketched in Algorithm \ref{alg_red} \footnote{A Matlab implementation of both the full and reduced matrix solvers for the discrete NS equation can be downloaded from https://sites.google.com/view/gerhard-kirsten/software upon acceptance of this article.}.

\begin{algorithm}[H]
\caption{Reduced model Algorithm}
\begin{algorithmic}[1]
\State Consider the projected initial conditions $(\widehat{{\bm U}}(0)$, $\widehat{{\bm V}}(0))$ and choose the number of time steps $n_t$  
\For{$i=0,\ldots, n_t-1$}
\State Solve \eqref{pres_corr_red} for $\widehat{\bm P}^{(i+1)}$ \label{step1}
\State Evaluate $\widehat{F}_U(\widetilde{\bm U},\widetilde{\bm V},t)$ and $\widehat{F}_V(\widetilde{\bm U},\widetilde{\bm V},t)$ with DEIM \label{step2}
\State Solve \eqref{nstimedisc} for $\widehat{\bm U}^{(i+1)}$ and $\widehat{\bm V}^{(i+1)}$ \label{step3}
\EndFor
\end{algorithmic}
\label{alg_red}
\end{algorithm}

\section{The tree structure algorithm for the NS equation}\label{sec:tsans}

In this section we are going to couple the multilinear approximation of the Navier-Stokes equation with the Tree Structure Algorithm (TSA), an algorithm to approximate the HJB equation arising from the optimal control problem. We will first introduce briefly the general procedure for the TSA and next we are going to present the coupling of these two techniques.



\subsection{Dynamic programming on a tree structure}\label{ssec:tsans}
We introduce the essential ingredients of the DP approach based on a tree built on the discrete dynamical system. The interested reader will find more details on the topic in \cite{afs19}.\\

We consider the discrete approximation of the DP principle. Fixed the number of time steps $n_t$ and the time step $\Delta t: = [(T-t)/n_t]$, the discrete DP reads
\begin{equation}\label{SL}
\left\{\begin{array}{ll}
V^{n}=\min\limits_{a\in A}[\Delta t\,L(x, a, t_n)+e^{-\lambda \Delta t}V^{n+1}(x+\Delta t f(x, a, t_n))], \quad
 n= n_t-1,\dots, 0,\\
V^{n_t}=g(x),\hspace{7.5cm} x \in \R^d,
\end{array}\right.
\end{equation}
where $t_n=t+n \Delta t,\, t_{\overline N} = T$, and $V^n:=V(x, t_n).$
We discretize the control set $A$ with step-size $\Delta a$ obtaining a discrete control set with a finite number of controls $A^{\Delta u}=\{a_1,...,a_M \}$. In what follows we denote by $A$ the discrete set to ease the notation.

Now, we start from the initial condition $x$ and we follow the  discrete dynamics employing the explicit Euler scheme and $M$ discrete controls $a_j$ 
\begin{equation}\label{def:ddyn}
\zeta_j^1 = x+ \Delta t \, f(x,a_j,t_0),\qquad j=1,\ldots,M.
\end{equation}

Therefore, denoting  the root of the tree with $\mathcal{T}^0=\{x\}$, we get the first level of the tree  $\mathcal{T}^1 =\{\zeta_1^1,\ldots, \zeta^1_M\}$. The procedure can be iterated so that the $n-$th {\em time level} will be given by 
$$\mathcal{T}^n = \{ \zeta^{n-1}_i + \Delta t f(\zeta^{n-1}_i, a_j,t_{n-1}) \}_{j=1}^{M}\quad i = 1,\ldots, M^{n-1}.$$ 
and the entire tree can be represented as
 $$\mathcal{T}:= \{ \zeta_j^n  \}_{j=1}^{M^n},\quad n=0,\ldots n_t,$$ 
where $\zeta^n_i$ is the evolution of the dynamics at time $t_n$ using the controls $\{a_{j_k}\}_{k=0}^{n-1}$:
$$\zeta_{i_n}^n = \zeta_{i_{n-1}}^{n-1} + \Delta t f(\zeta_{i_{n-1}}^{n-1}, a_{j_{n-1}},t_{n-1})= x+ \Delta t \sum_{k=0}^{n-1} f(\zeta^k_{i_k}, a_{j_k},t_k), $$
with $\zeta^0 = x$, $i_k = \floor*{\dfrac{i_{k+1}}{M}}$ and $j_k\equiv i_{k+1} \mbox{mod } M$.

Although the TSA allows to deal with high dimensional problems, the cardinality of tree grows exponentially in the time steps and in the number of nodes, $i.e.$ $|\mathcal{T}|=O(M^{n_t})$, yielding problems in the memory allocations.
For this reason we introduce a \emph{pruning criteria} based on the distance between nodes.
Therefore, two nodes $\zeta^n_i$ and $\zeta^n_j$ will be merged if 
\begin{equation}\label{tol_cri}
\Vert \zeta^n_i-\zeta^n_j \Vert \le \ep, \quad \mbox{ with }i\ne j \mbox{ and } n = 0,\ldots, n_t, 
\end{equation}
for a given threshold $\ep>0$. In \cite{saf21} the authors show that the threshold $\ep>0$ must scale quadratically in the time steps to ensure first order convergence. 

Once constructed the tree $\mathcal{T}$, we can pass to the computation of the numerical value function $V(x,t)$. The TSA defines a time dependent grid $\mathcal{T}^n=\{\zeta^n_j\}_{j=1}^{M^n}$ for $n=0,\ldots, n_t$ and
\eqref{HJB} can be approximated as follows: 
\begin{equation}
\begin{cases}
V^{n}(\zeta^n_i)= \min\limits_{a\in A} \{e^{-\lambda \Delta t} V^{n+1}(\zeta^n_i+\Delta t f(\zeta^n_i,a,t_n)) +\Delta t \, L(\zeta^n_i,a,t_n) \}, \\
\qquad\qquad \qquad\qquad \qquad\qquad \qquad\qquad \qquad\qquad  \zeta^n_i \in \mathcal{T}^n\,, n = n_t-1,\ldots, 0, \\
V^{n_t}(\zeta^{n_t}_i)= g(\zeta_i^{n_t}), \qquad\qquad \qquad\qquad \qquad\qquad   \zeta_i^{n_t} \in \mathcal{T}^{n_t}.
\end{cases}
\label{HJBt2}
\end{equation}

The minimization in \eqref{HJBt2} is solved by comparison on the discrete set $A$. 

\subsection{Coupling TSA and 2S-POD-DEIM}\label{ssec:test1}

Introduced the main ingredients for the multinear approximation of the NS equation and the TSA, in this section we are going to show how to couple these two techniques in order to solve the optimal control problem. The procedure is divided into two steps: an offline and an online phase.
\begin{itemize}
    \item {\bf Offline Phase}
    
    In the offline phase we build the 2S-POD-DEIM basis and construct the reduced dynamics which will be employed in the online phase. In this step we explore the manifold of possible evolutions of the controlled dynamics and we aim to capture the main features of the dynamical system. We first fix a time step $\widehat{\Delta t}$ and $\widehat{M}$ number of discrete controls. 
    Following Algorithm \ref{alg_full}, the tree structure is constructed in the full dimension with the fixed parameters $\widehat{\Delta t}$ and $\widehat{M}$. Since at this stage the problem is high-dimensional, few time steps and few controls will be selected for the offline phase. Once the tree has been built, we can pass to the construction of the 2S-POD-DEIM basis. Following Section \ref{2spod}, we introduce a truncating tolerance $tol$ and we generate the reduced basis $\bm{U}_\ast \in \RR^{n_x - 1 \times k_{1,\ast}}$, $\bm{V}_\ast \in \RR^{n_y - 1 \times k_{2,\ast}}$, $\bm{P}_\ast \in \RR^{n_x \times k_{3,\ast}}$,$\bm{\Phi}_{\bm{U}, \ast} \in \mathbb{R}^{n_x -1 \times p_{1,\ast}}$, $\bm{\Phi}_{\bm{V}, \ast} \in \mathbb{R}^{n_y -1 \times p_{2,\ast}}$, and the permutation basis $\bm{D}_{\bm{U}, \ast} \in \mathbb{R}^{n_x -1 \times p_{1,\ast}}$  and $\bm{D}_{\bm{V}, \ast} \in \mathbb{R}^{n_y -1 \times p_{2,\ast}}$, where $\ast = \{\ell,r\}$. Furthermore, all low-dimensional time-independent matrices are constructed and stored in this phase.

    \item {\bf Online Phase}
    
    In the online phase we solve the optimal control problem via TSA directly on the reduced model constructed in the offline phase. Since we reduced the dimension of the system, in this step it is possible to consider more time steps and/or more discrete controls with respect to the offline phase. First of all we construct the reduced tree following Algorithm \ref{alg_red}. The nodes of the tree represent the grid for the numerical resolution of the DDP \eqref{HJBt2}. At this point we can solve the DPP on the tree structure, obtaining the discrete value function $\{V_{red}^n(\zeta_i^n)\}_{i,n}$ on the tree. The last part of this phase concerns the reconstruction of the control signal and the controlled trajectory.
    Starting from the initial condition, $i.e.$ the root of the tree, we can follow the branches returning the minimum

    \begin{equation*} \label{feed:tree-POD}
\alpha^{n}_{*}:=\argmin\limits_{a\in U} \left\{ e^{-\lambda \Delta t}V_{red}^{n+1,\ell}(\zeta_{red}^{n}+\Delta t f_{red}(\zeta_{red}^{n},a,t_n)) +\Delta t \, L_{red}(\zeta_{red}^{n},a,t_n) \right\}.
\end{equation*}

    The computed control signal can now be plugged into the full dimensional dynamics to obtain the optimal trajectory in the original dimension.
\end{itemize}

\section{Numerical experiments}\label{sec:nens}

In this section we are going to test the proposed technique in different settings. In the first test we compare the performances between the full dimension model and the low dimension one in terms of CPU time. Moreover, the vector and the matricial cases will be compared. In the second and third tests we pass to the optimal control problem of the NS equation, in which we are interested in reaching a particular solution target: the stationary configuration. More precisely, in the second test the control will act on the entire domain via the use of a shape function, while in the third example the control operates on a subdomain $w$ located at the center of the domain. In the last example the control will operate on the boundaries and the reference solution will be represented by the trajectory obtained using a prefixed control. For all the numerical tests we will fix the space domain $\Omega = [0,1]^2$ and Reynolds number $r=100$. We also give some indications on the efficiency of the pruning technique via the {\em Pruning Ratio ($Ratio_p$)} that is defined as the ratio between the cardinality of the full tree and the cardinality of the pruned tree, $i.e.$
$$
Ratio_p = \frac{(M^{n_t+1}-1)/(M-1)}{|\mathcal{T}_M|}
$$
where $\mathcal{T}_M$ is the tree constructed using $M$ discrete controls under the pruning criteria.

\subsection{Test 1: Comparison full/low dimension}\label{ssec:test1}
In this example we investigate the efficiency of, not only discretizing the NS equation in matrix form, but also reducing the dimension of the resulting matrix equation by 2S-POD-DEIM.
 To this end we consider the NS equation \eqref{nscont}, and we fix $T=20$, $\Delta t=0.05$ and the tolerance for 2S-POD-DEIM equal to $10^{-3}$. We consider as initial condition $u_0 = v_0 \equiv 0$ and the following boundary conditions
$$
u (t,x)=v(t,x)=1, \quad (x_1,x_2)\in [0,1] \times \{1\}, t \in (0,T]
$$
and homogeneous Dirichlet conditions on the other walls.

\begin{figure}[htb!]
\centering
\minipage{0.4\textwidth}
  \includegraphics[width=\linewidth]{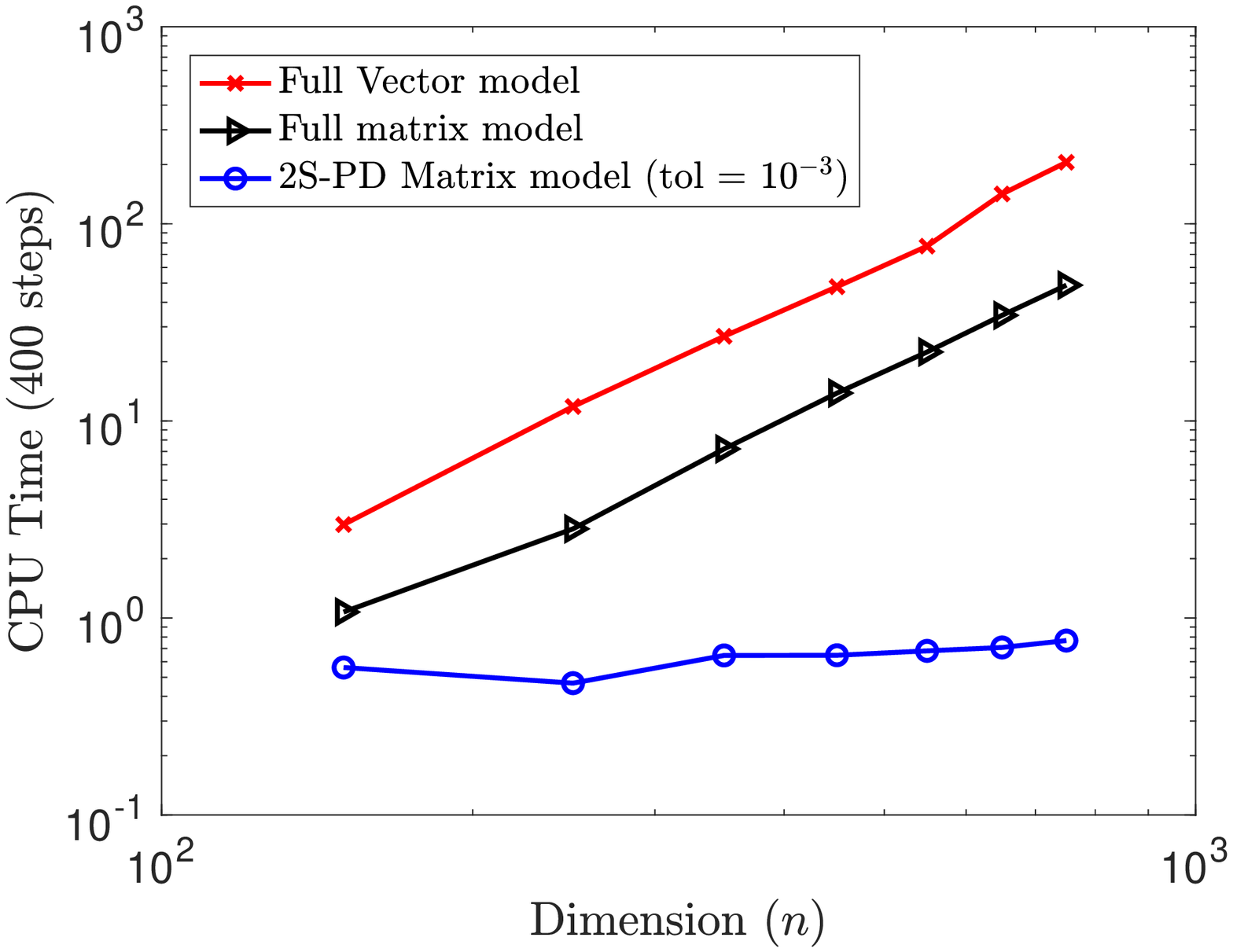}
\endminipage
\minipage{0.4\textwidth}
  \includegraphics[width=\linewidth]{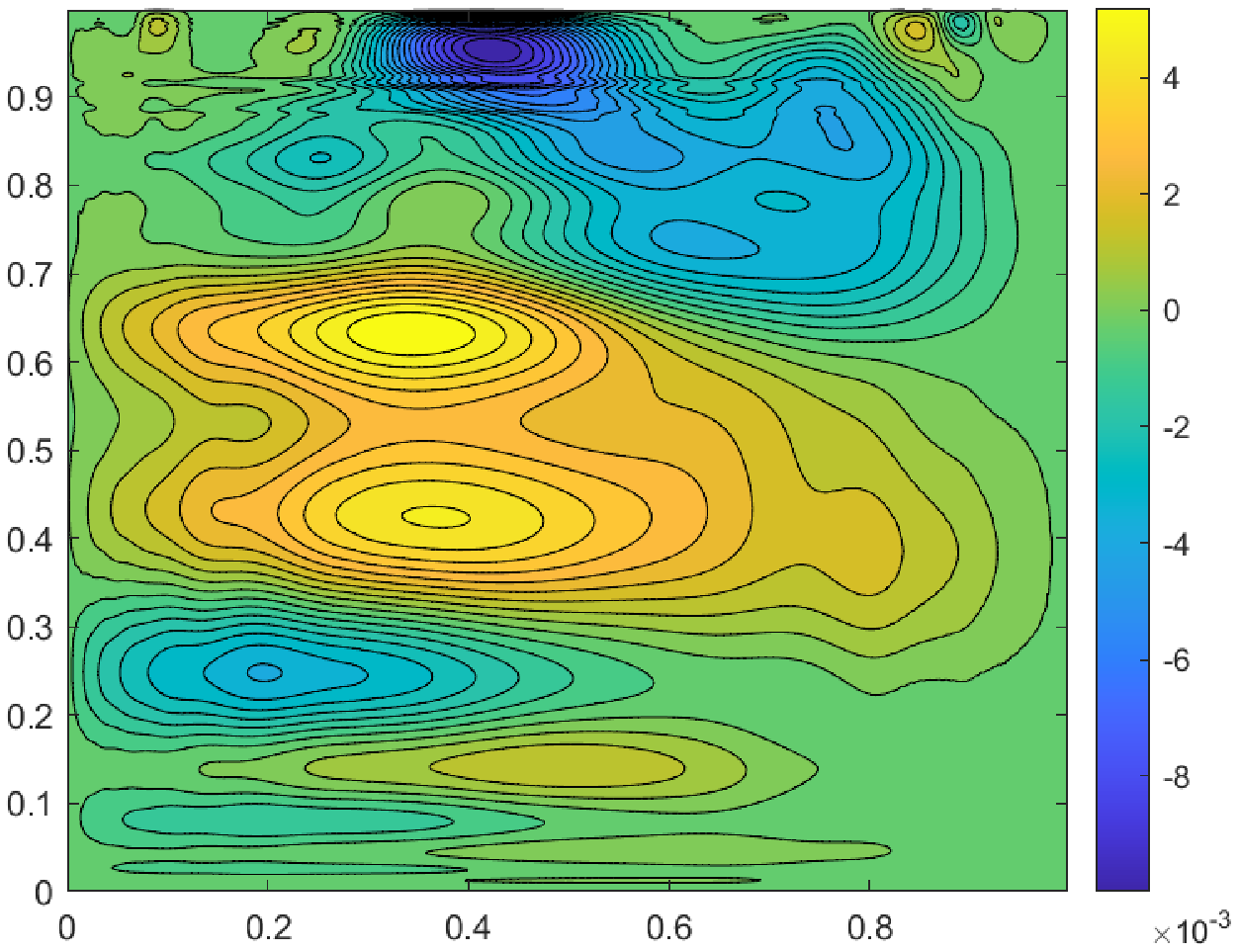}
\endminipage\hfill

\caption{Test 1: Comparison between the full vector model, the full matrix model and the 2S-POD-DEIM matrix model with $tol=10^{-3}$ in terms of CPU time varying the space dimension (left) and difference between the full solution and the lifted solution at final time with $n=750$ (right).}	
\label{full_red}
\end{figure}

For the experimental setup we consider $n = n_x = n_y$, for $n \in \{ 150, 250, 350, 450, 550, 650, 750\}$ and measure the CPU time required to evaluate the discrete NS equation at $n_t = \frac{T}{ \Delta t} = 400$ timesteps for the full dimensional vector model, the full dimensional matrix model and the 2S-POD-DEIM reduced matrix model. The results for the full dimensional vector model are obtained by running the Matlab software from \cite{seibold08}. The results are plotted in Figure \ref{full_red}, left.

It can be deduced from the plot that solving the full dimensional discrete NS equation in matrix form as opposed to vector form results in a good computational gain. This behaviour is typical due to the efficiency of the a-priori eigenvalue decomposition resulting in a simple solve by substitution for the Sylvester equations at each time step; see, e.g., \cite{dautilia20, kirsten_phd}. Furthermore, as expected, we notice that the reduced model is several orders of magnitude faster than both full-order model, and the nearly-constant timings for increasing $n$ indicates that the computational cost for the reduced model is indeed completely independent of the full dimension $n$. 

 In the right panel of Figure \ref{full_red} we present the difference 
$$
{\bm U}(T)-{\bm U}_\ell \widehat{{\bm U}}(T) {\bm U}^\top_r
$$
which describes the error related to the projection at the final time, having fixed $n=750$. It is clear from the plot that the projection error has the same order of the chosen tolerance, $i.e.$ almost $10^{-3}$. At the top wall the projection error is slightly higher due to the non-homogeneous boundary condition.

These promising results show that a very fine discretization of the NS equation can be solved in a fraction of a second online with a small projection error of order $10^{-3}$.

\subsection{Test 2: The problem of long time behaviour of the solution}\label{ssec:test2}

\begin{figure}[htb!]
\minipage{0.33\textwidth}
  \includegraphics[width=\linewidth]{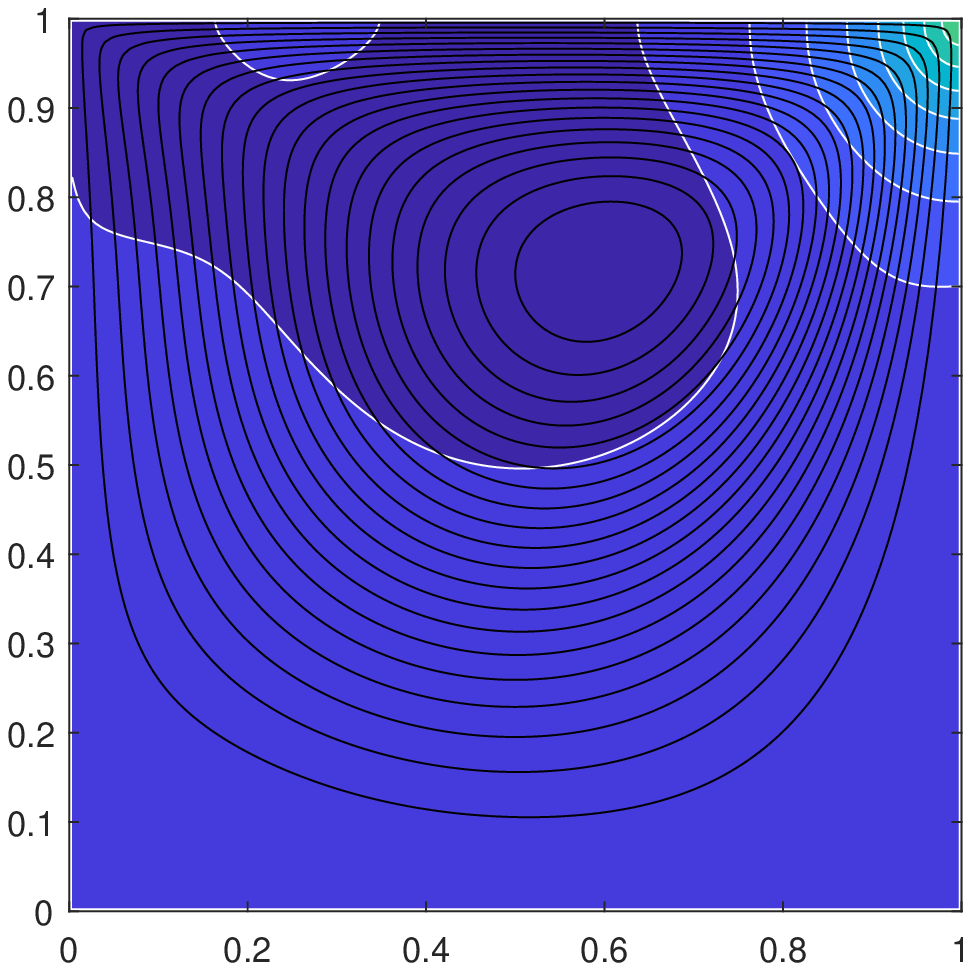}
\endminipage\hfill
\minipage{0.33\textwidth}
  \includegraphics[width=\linewidth]{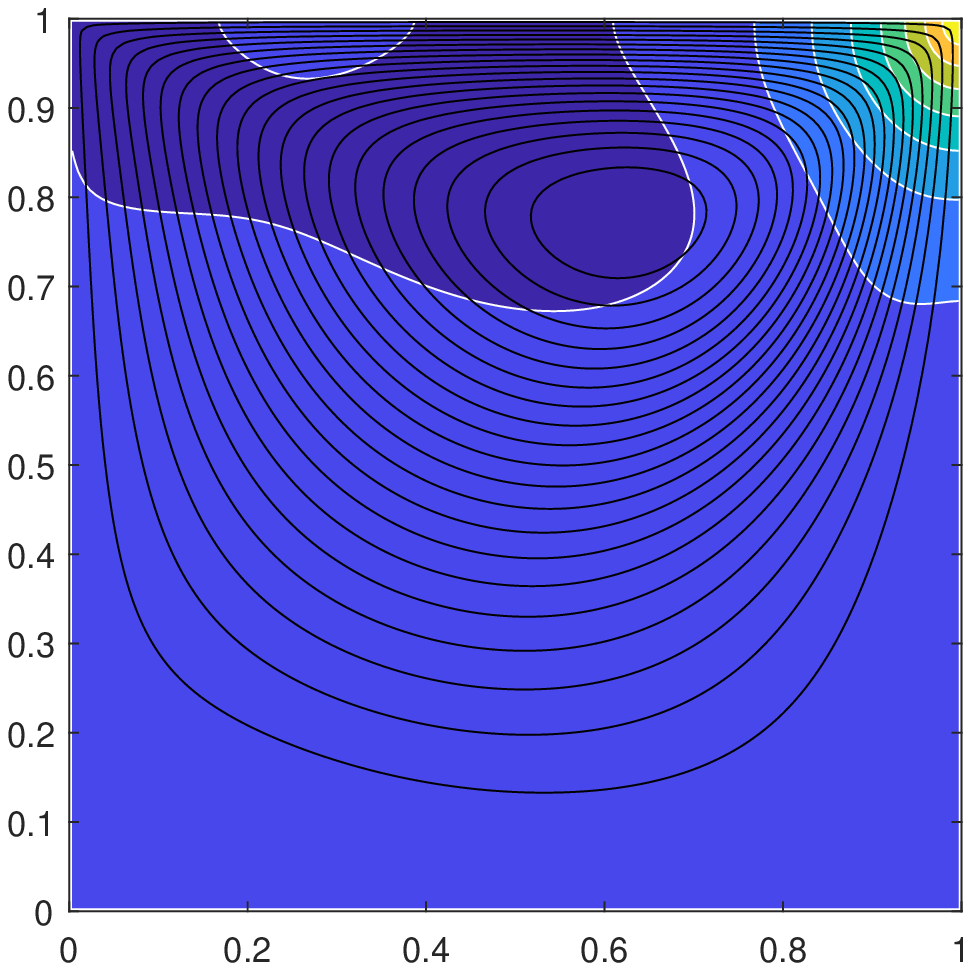}
\endminipage\hfill
\minipage{0.33\textwidth}
  \includegraphics[width=\linewidth]{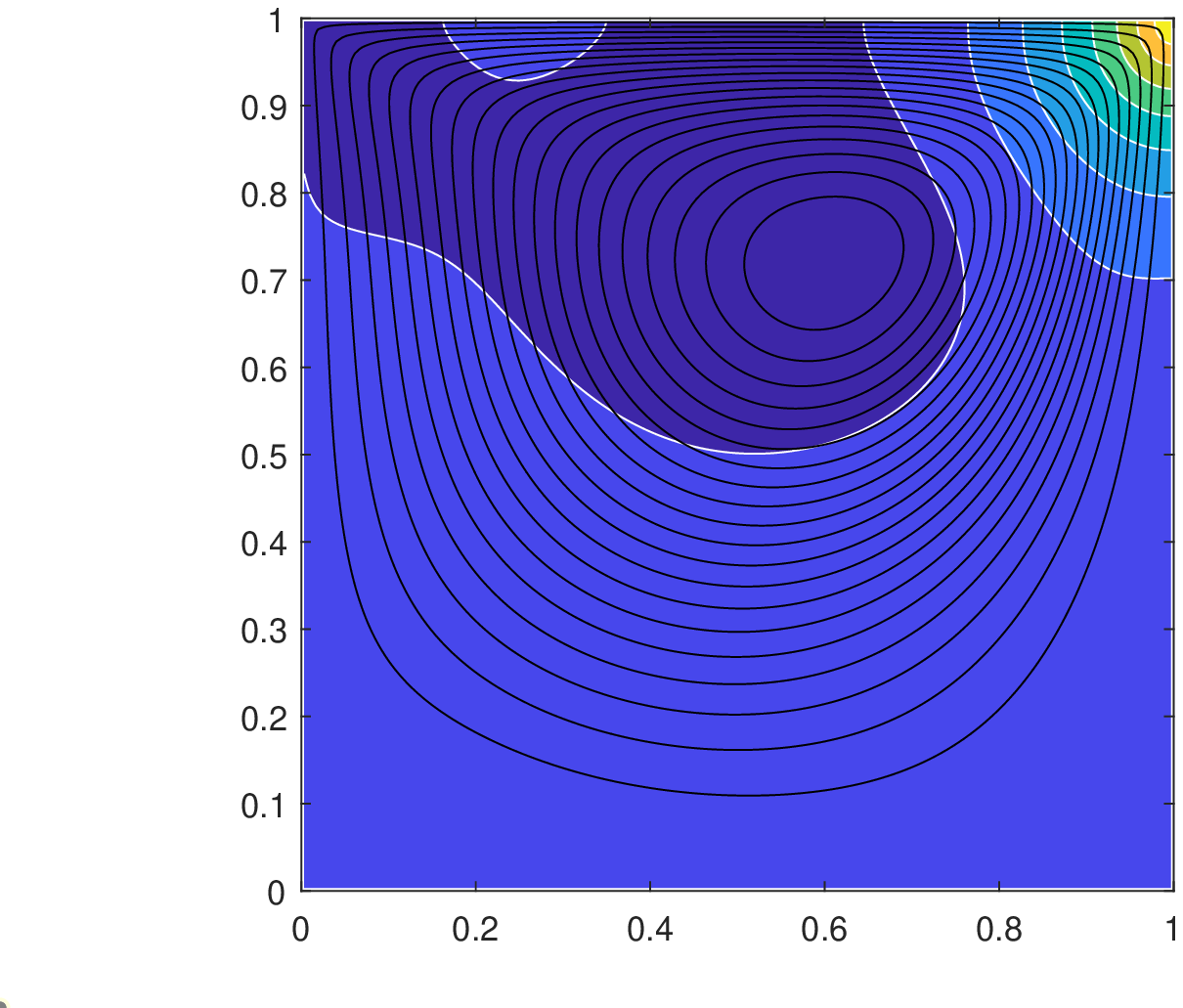}
\endminipage
\caption{Test 2: Stationary solution ($t=20$) (left),  uncontrolled solution at $t=2$ (central) and controlled solution at $t=2$ (right) for 2S-POD and $n=201$.}	
\label{figtransport}
\end{figure}

\begin{figure}[htb!]
\minipage{0.33\textwidth}
  \includegraphics[width=\linewidth]{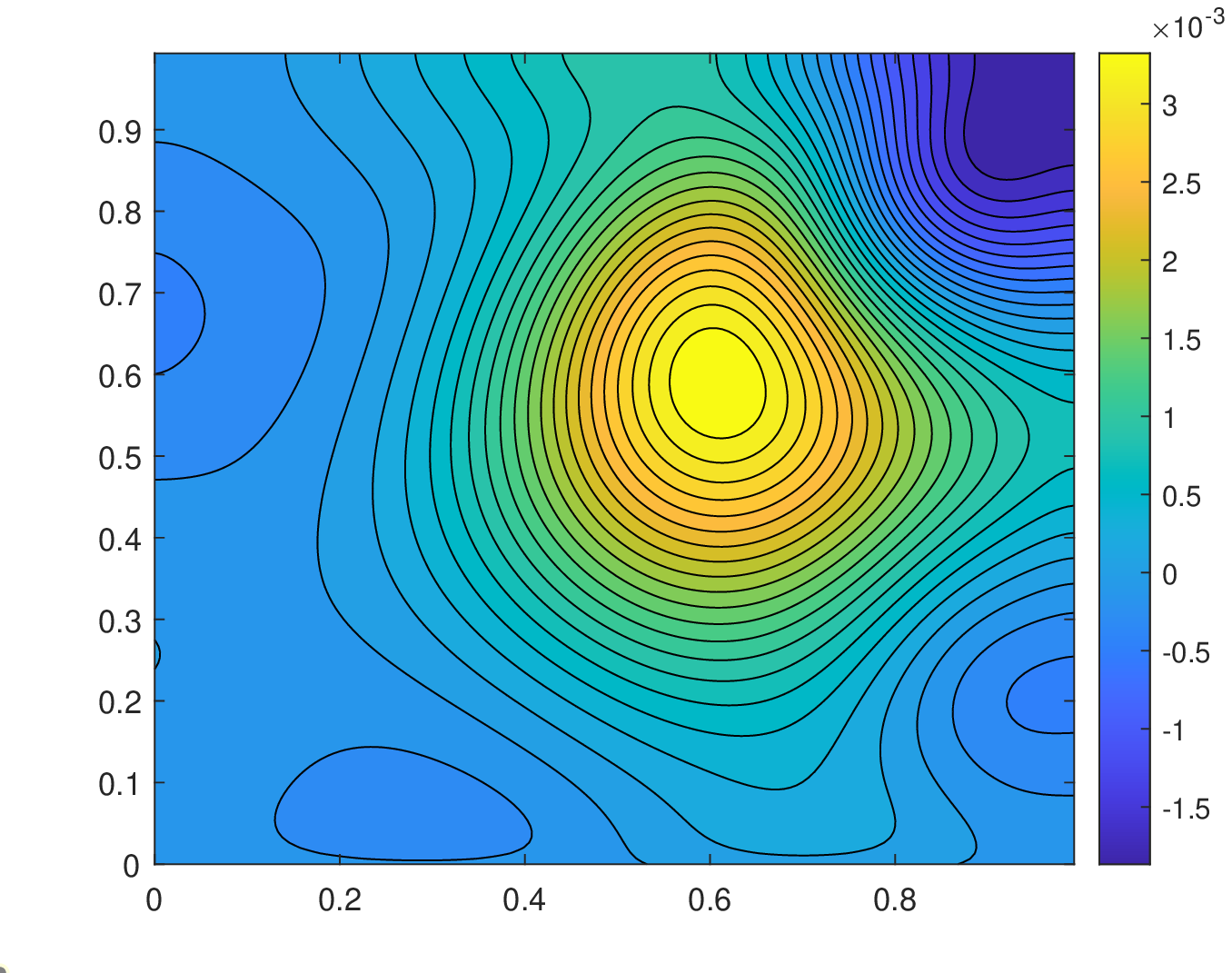}
\endminipage\hfill
\minipage{0.33\textwidth}
  \includegraphics[width=\linewidth]{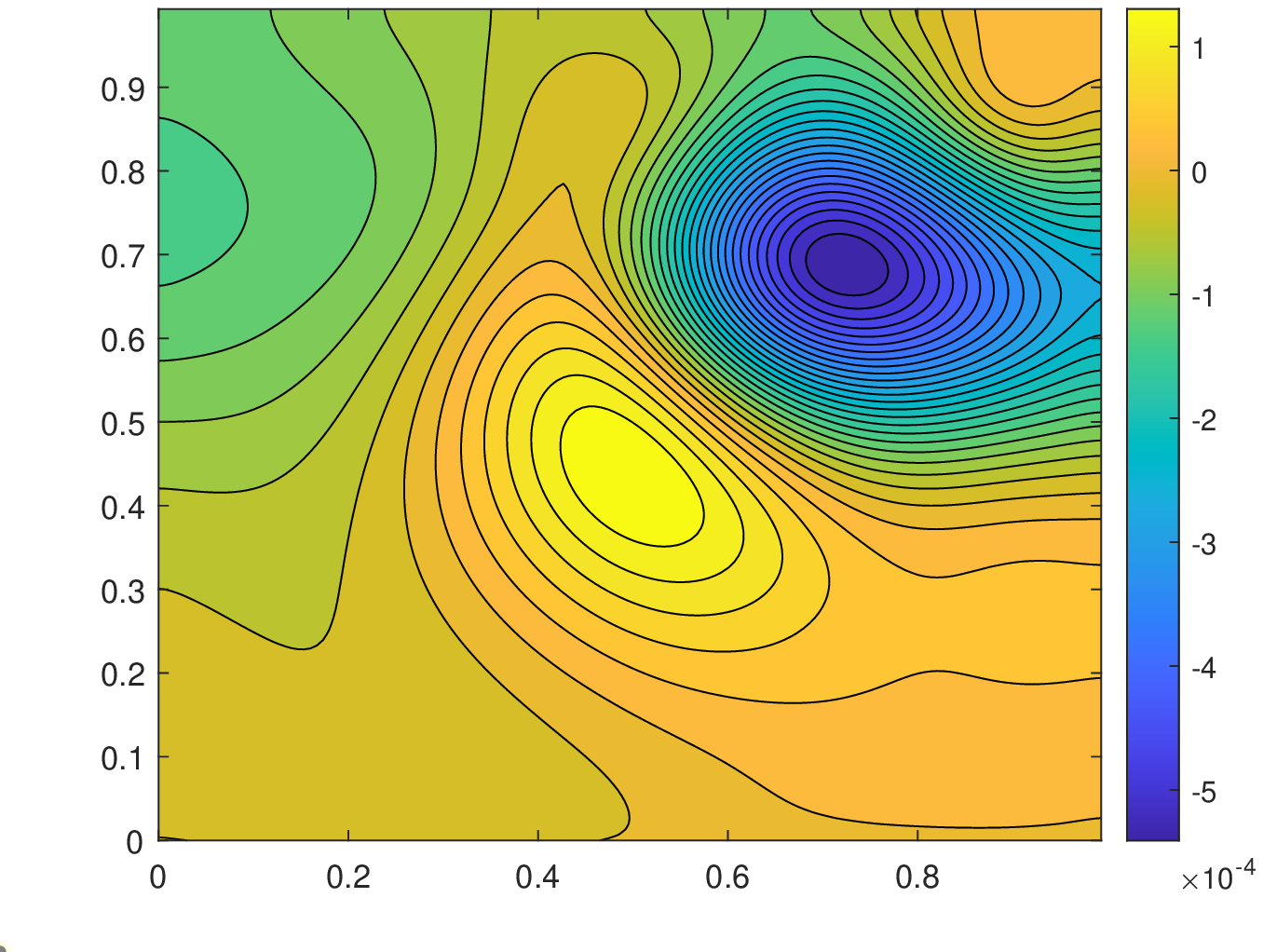}
\endminipage\hfill
\minipage{0.33\textwidth}
  \includegraphics[width=\linewidth]{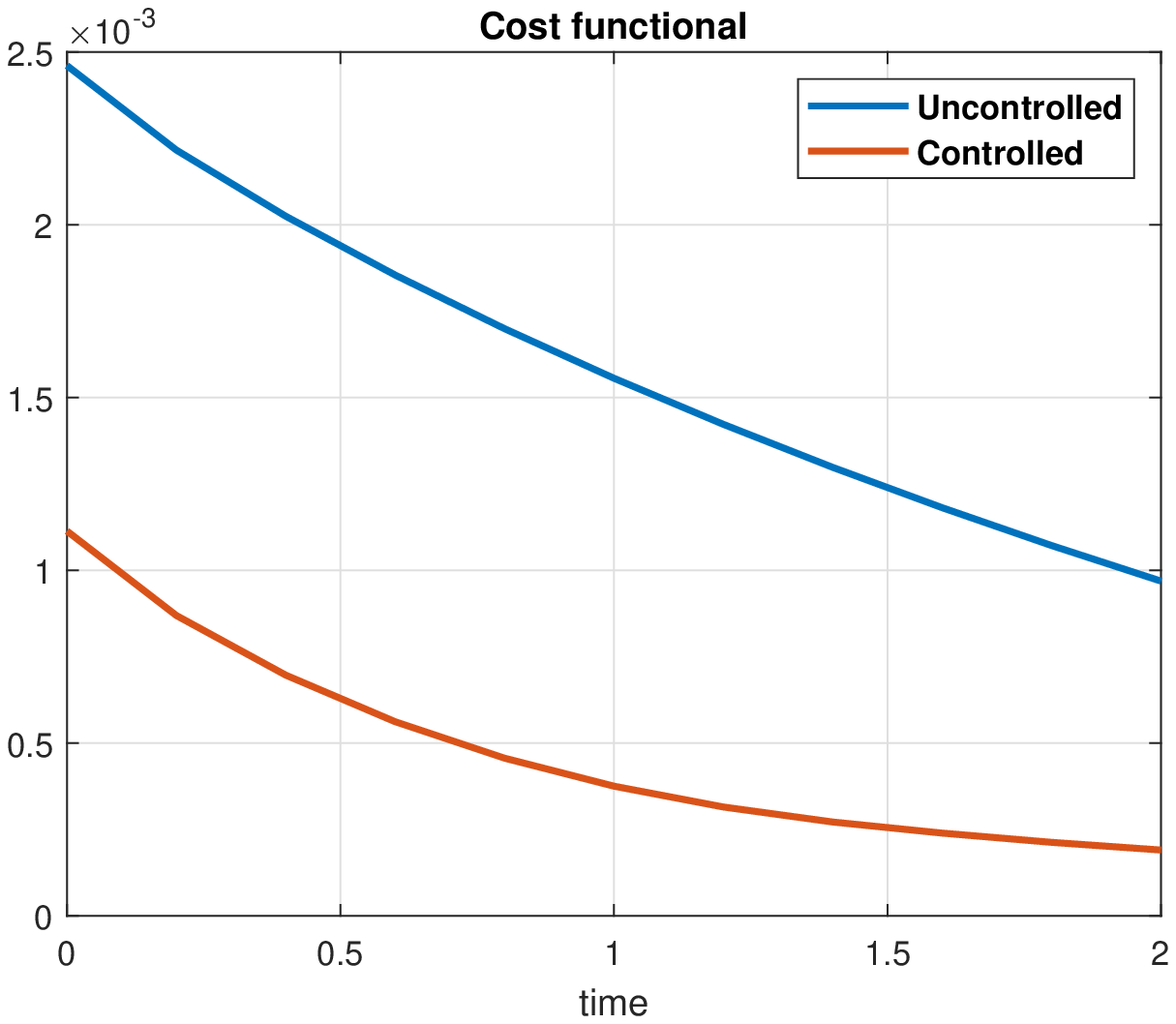}
\endminipage\hfill
\caption{Test 2: Pressure difference between stationary and uncontrolled solution (left), between stationary and controlled solution (central) at the final time and cost functional (right) for 2S-POD and $n_x=201$.}	
\label{fig_test2_3}
\end{figure}

In this test we want to reach a target solution acting on a scalar control that appears in the Navier-Stokes equation as an additional term as in \eqref{eq:controlpb1}. The initial condition and the boundary conditions coincide with the previous example. We are interested in the stationary solution of the Navier-Stokes equation and it can be obtained analysing the long time behaviour of the uncontrolled dynamics. In this case the stationary solution $\overline{y}$ has been fixed as the uncontrolled solution at time $t=20$, since the solution does not present relevant changes for larger time intervals. The optimal control problem is based on the minimization of the following cost functional
$$
J(\alpha)= \int_0^T \| y(\cdot, t;\alpha) -\overline y(\cdot, t)\|^2_{L^2(\Omega)} + \gamma \|\alpha\|^2) e^{-\lambda t} dt +\| y(\cdot, T;\alpha) -\overline y(\cdot, T)\|^2_{L^2(\Omega)}
$$
with $\gamma = 10^{-3}$.
The parameters of the discrete problem are the following: $\Delta t=0.1$, $A=\{0,0.5,1\}$ and $T=2$. A similar example has been studied in \cite{AH2014}. The geometric pruning criteria will be applied by selecting the threshold $\epsilon_{\mathcal{T}}= \Delta t^2$ to ensure first order convergence. We fix the gridpoints $n_x=n_y=201$ and we apply the 2S-POD-DEIM technique with tolerance $\tau=10^{-3}$. The reduction techniques provide the following basis: $\bm{U}_\ell \in \mathbb{R}^{200 \times 61}$, $\bm{U}_r \in \mathbb{R}^{200 \times 57}$, $\bm{V}_\ell, \bm{V}_r, \bm{P}_{ \ell}, \bm{P}_{ r}  \in \mathbb{R}^{200 \times 70}$, $\bm{\Phi}_{\bm{U}, \ell}, \bm{\Phi}_{\bm{V}, \ell} \in \mathbb{R}^{200 \times 62}$ and $\bm{\Phi}_{\bm{U}, r},\bm{\Phi}_{\bm{V}, r} \in \mathbb{R}^{200 \times 52}$. In Figure \ref{figtransport} we display respectively the pressure field computed in the stationary case (left panel) and the uncontrolled (central panel) and the controlled solution (right panel) at time $t=2$. We show the contour lines of the pressure field and the closed contour lines of the stream function. It is possible to notice visually how the controlled dynamics looks similar to the stationary solution, while the uncontrolled is still far from the asymptotic behaviour. The difference in the pressure field at the final time between the uncontrolled solution and $\widetilde{y}$ is displayed in the left panel of Figure \ref{fig_test2_3}, while in the central panel we show the difference between the controlled and stationary solution. We notice that the order in latter case is $\approx 10^{-4}$, while in the first case is $\approx 10^{-3}$, demonstrating how the solution of the optimal control converges more rapidly to the stationary configuration. The right panel of Figure \ref{fig_test2_3} shows the comparison of the cost functional in the controlled and uncontrolled setting, where we can see again the faster convergence of the controlled dynamics to the stationary solution. In this case the application of the pruning criteria yields a cardinality of the tree equal to $2994$, whereas the cardinality of the full tree is $88573$, corresponding to a pruning ratio of almost 30.

\

\subsection{Test 3: control on an internal subdomain $\omega\subset \Omega$}\label{ssec:test3}

In this experiment our aim is to reach a target solution acting on a scalar control that appears in the Navier-Stokes equation as an additional term concentrated on a subdomain $\omega\subset\Omega$ as in  \eqref{eq:controlpb2}. 
We consider $\omega=[0.3,0.7]^2$, $i.e.$ the control will operate on a central smaller square. In this example we consider homogeneous Dirichlet boundary conditions for all the walls and our scope is to drive the solution to the equilibrium $\overline{y} \equiv 0$. In this case we select a cost functional depending only on the final cost
$$
J(\alpha)= \| y(\cdot, T;\alpha)\|^2_{L^2(\Omega)}.
$$

\begin{figure}[htb!]
\minipage{0.33\textwidth}
  \includegraphics[width=\linewidth]{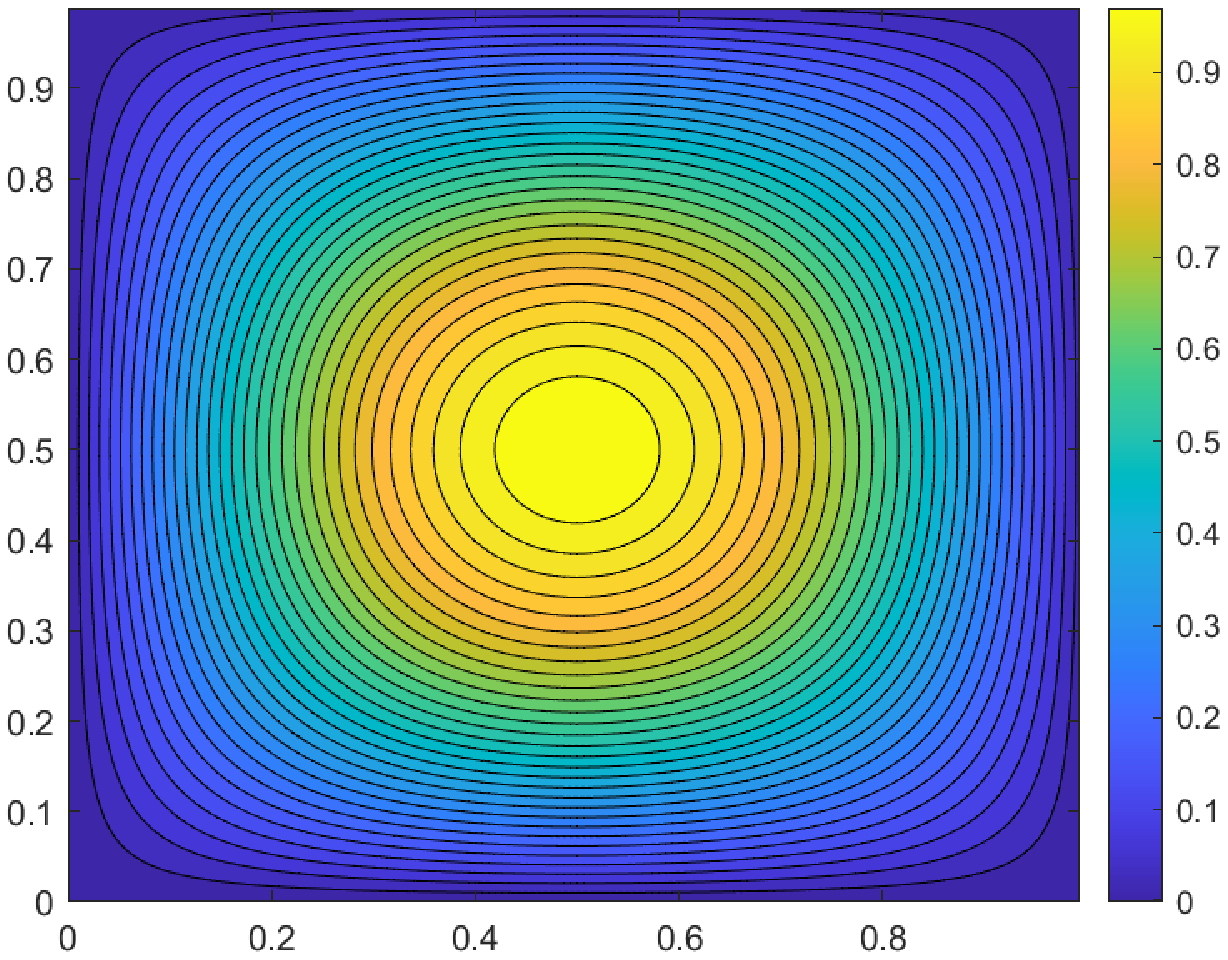}
\endminipage\hfill
\minipage{0.33\textwidth}
  \includegraphics[width=\linewidth]{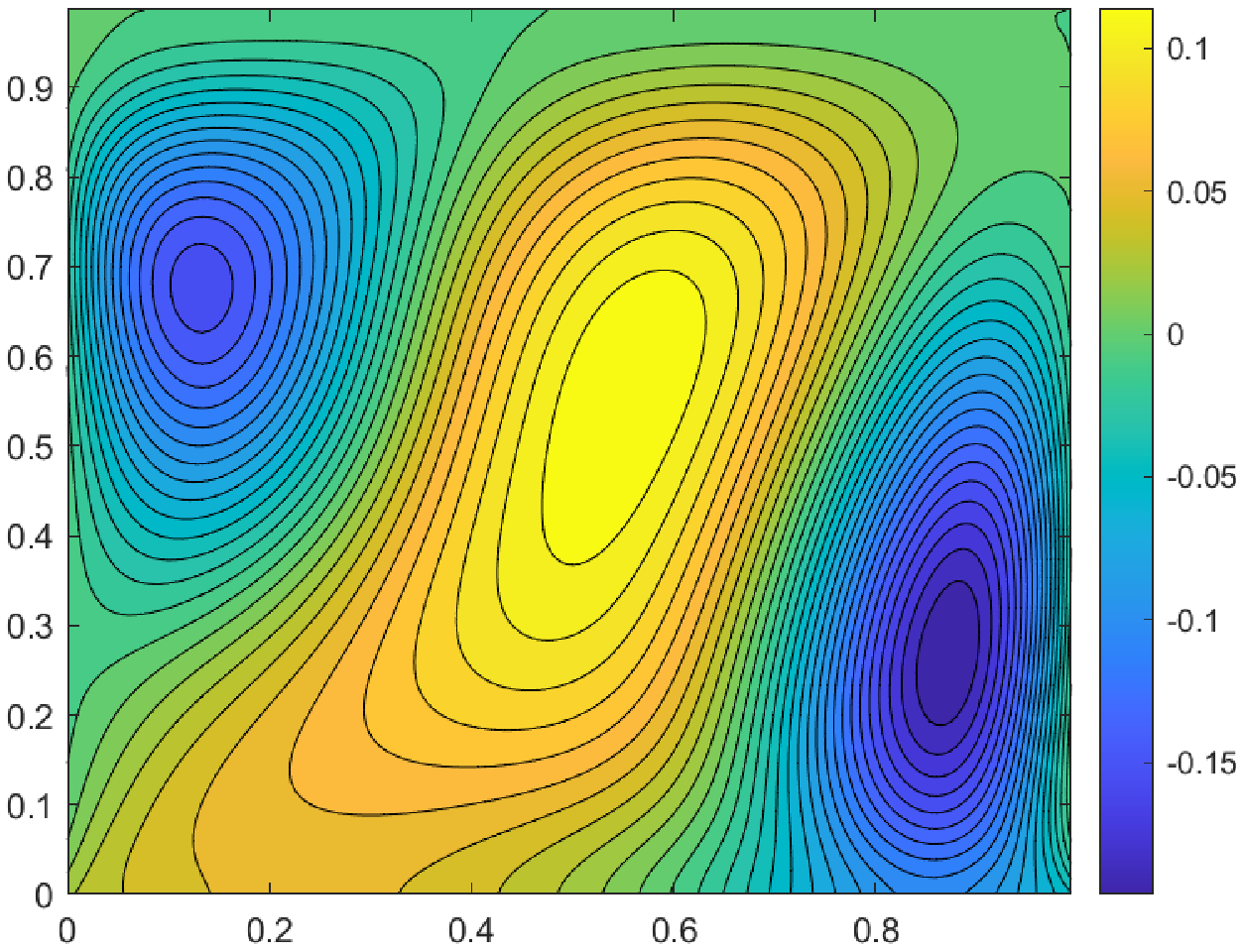}
\endminipage\hfill
\minipage{0.33\textwidth}
  \includegraphics[width=\linewidth]{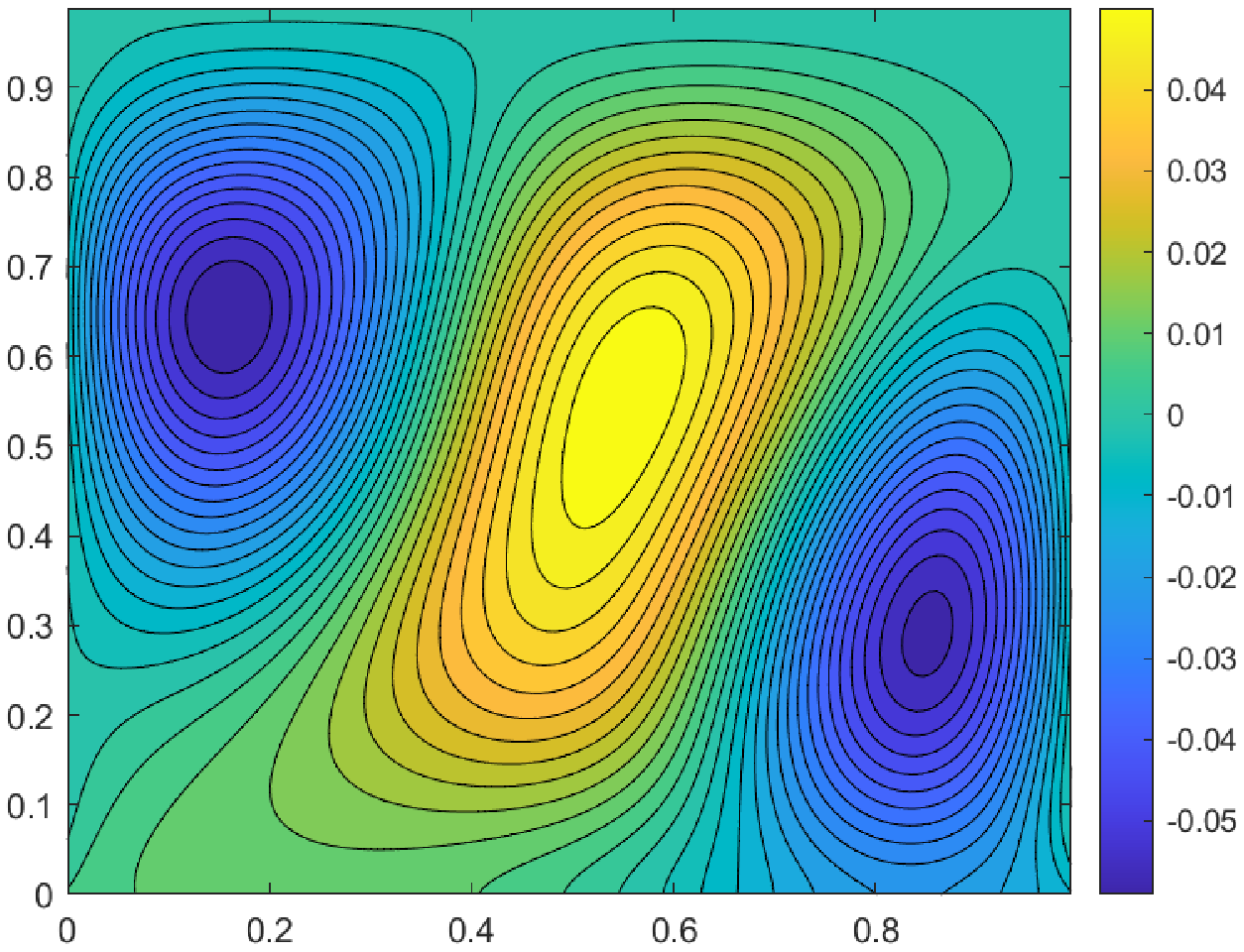}
\endminipage
\caption{Test 3: Uncontrolled solution at $t=0$ (left), $t=1$ (central) and $t=2$ (right) for 2S-POD and $n_x=201$.}	
\label{fig_test_2}
\end{figure}

We consider $\Delta t=0.1$, $T=2, A=[0,1]$ and we will vary the number of discrete controls. We consider the following initial condition
$$
u_0 = v_0 = \sin(\pi x) \sin(\pi y), \quad (x,y) \in [0,1]^2.
$$
In Figure \ref{fig_test_2} we show the behaviour of the uncontrolled solution ${\bm U}$ for different times. We note that the norm of the solution is decreasing due to the viscosity term. The aim of the corresponding optimal control problem is to accelerate this decay.

\begin{figure}[htb!]
\minipage{0.33\textwidth}
  \includegraphics[width=\linewidth]{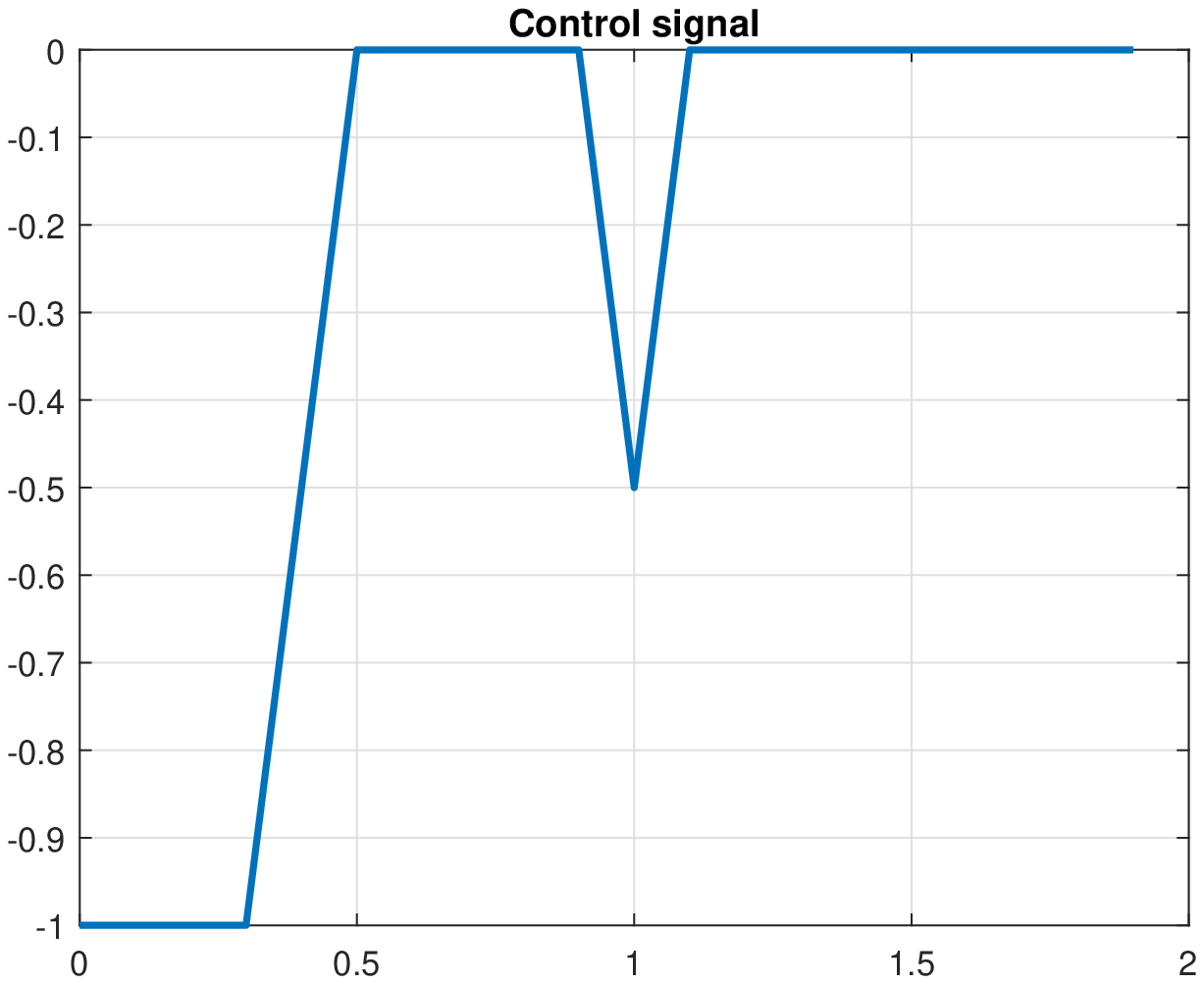}
\endminipage\hfill
\minipage{0.33\textwidth}
  \includegraphics[width=\linewidth]{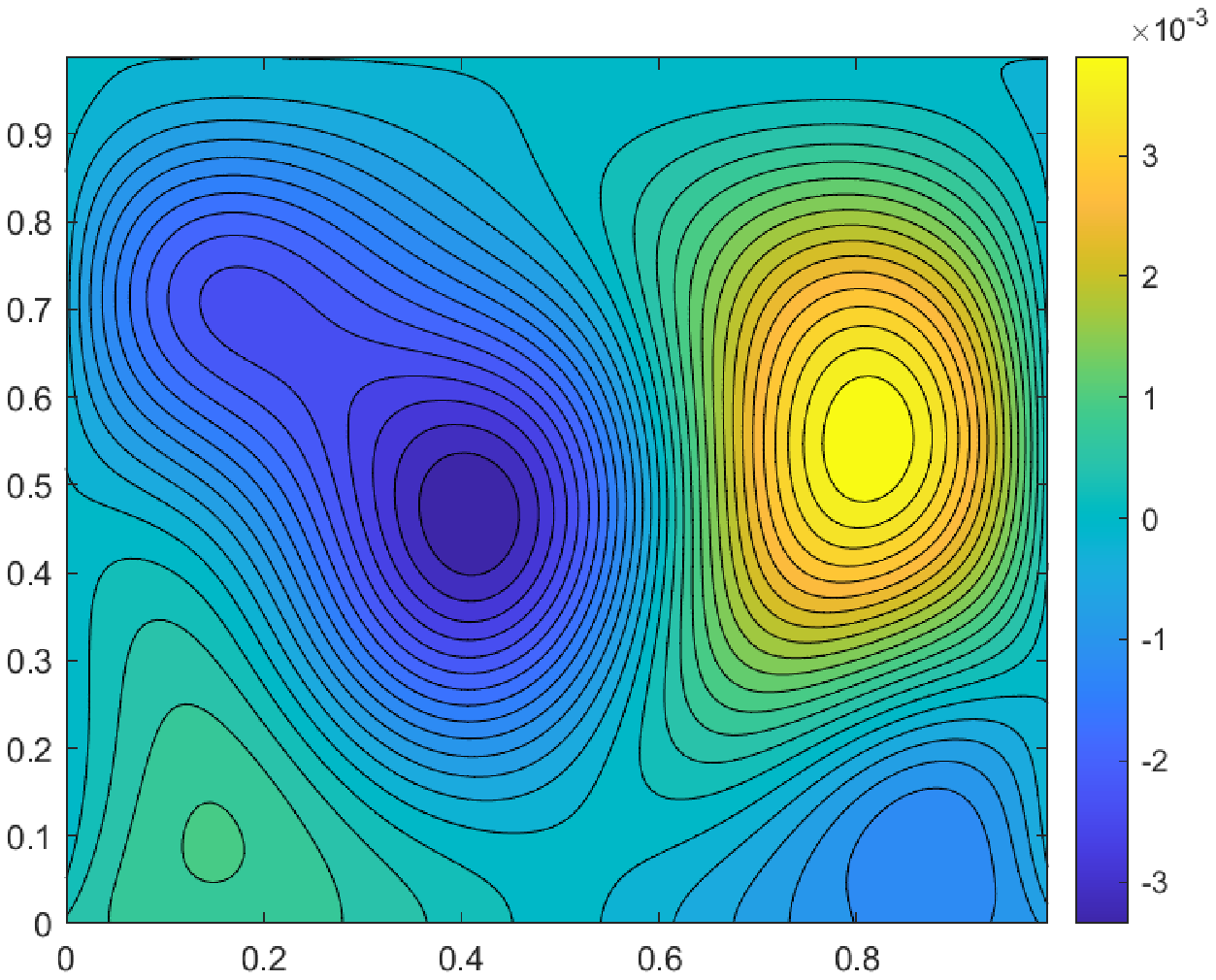}
\endminipage\hfill
\minipage{0.33\textwidth}
  \includegraphics[width=\linewidth]{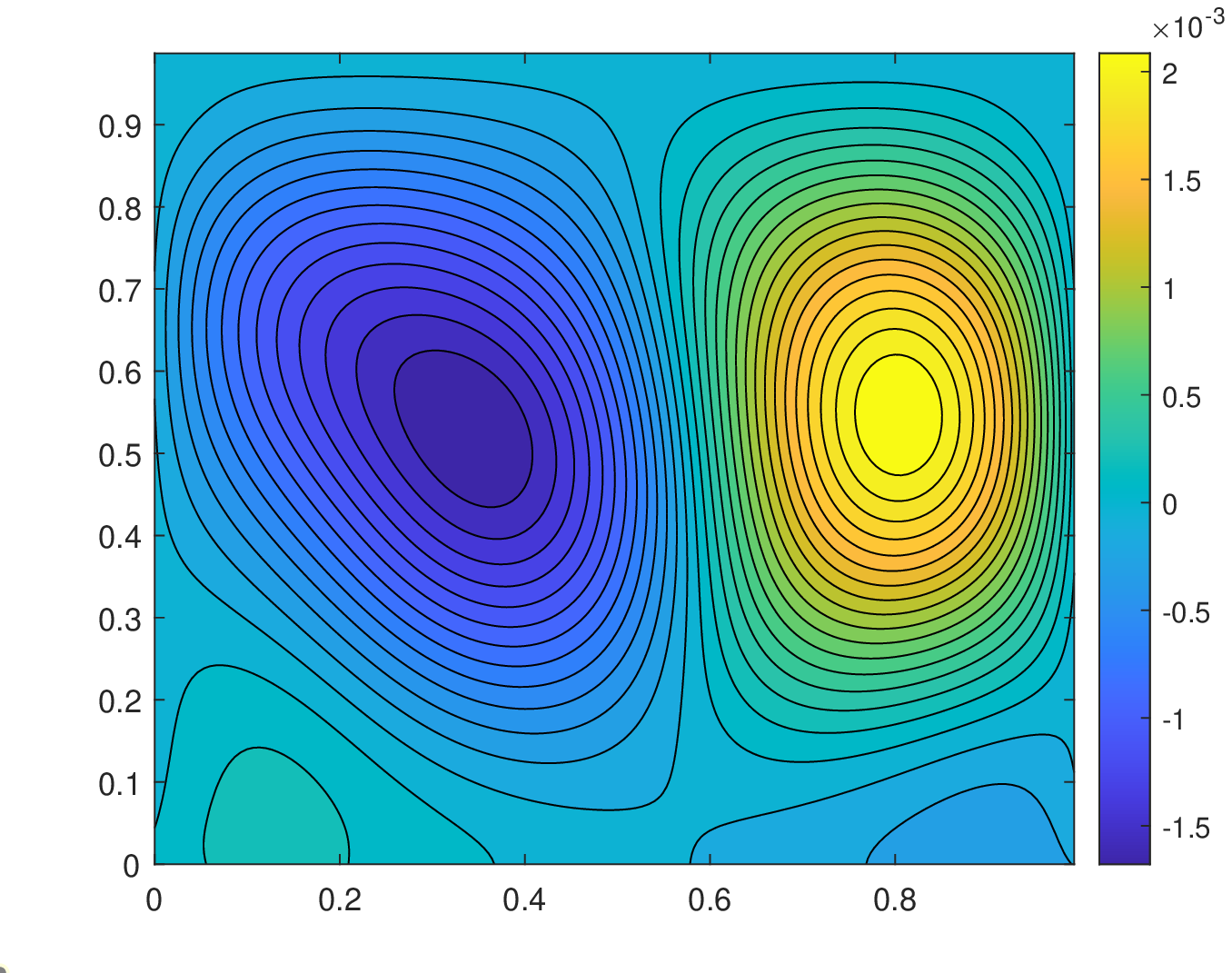}
\endminipage
\caption{Test 3: Control signal (left), controlled solution at $t=1$ (central) and $t=2$ (right) for 2S-POD and $n_x=201$.}	
\label{fig_test_22}
\end{figure}

We apply the 2S-POD-DEIM approach and we construct the following basis: $\bm{U}_\ell, \bm{U}_r, \bm{V}_\ell, \bm{V}_r  \in \mathbb{R}^{200 \times 68}$,  $\bm{P}_{ \ell},\bm{P}_{ r} \in \mathbb{R}^{200 \times 67}$ , $\bm{\Phi}_{\bm{U}, \ell}, \bm{\Phi}_{\bm{U},r}, \bm{\Phi}_{\bm{V}, \ell}, \bm{\Phi}_{\bm{V}, r} \in \mathbb{R}^{200 \times 70}$.
Figure \ref{fig_test_22} displays the results obtained by the coupling of the TSA and 2S-POD-DEIM. The left panel shows the control signal which presents a non-decreasing behaviour at the beginning of the time interval but starts oscillating in the middle. The central and the right panels we report the controlled solution at the time instances $t \in \{1,2\}$. We note that the maximum of the controlled solution is order $10^{-3}$ at time $t=1$, whereas for the uncontrolled dynamics it is stuck to $10^{-1}$. After $t=1$, the control stops acting and the decrease is just due to viscous term in the equation.
In Table \ref{Table_comparison} we present the comparison between the uncontrolled dynamics and the controlled solution varying the number of discrete controls. In term of the cost functional, the TSA gets almost one order of magnitude with respect to the uncontrolled case and we see an improvement increasing the number of controls. Moreover we report the cardinality of the tree coupled with the pruning technique. Looking at the P-Ratio we note the pruning criteria yields to a great benefit in terms of memory storage and this improvement increases as we consider more discrete controls.

	\begin{table}[htbp]							
		\centering
		\begin{tabular}{c |c c c}
					\toprule
					& Cost functional & Nodes & $Ratio_p$
					     \\
					\midrule
					Uncontrolled & 3.60e-4 & & \\
					Controlled $M=2$ & 5.44e-5 & 6064 & 345 \\
					Controlled $M=3$ &   5.27e-5 & 83273 & 6.2e4 \\
			\bottomrule		
		\end{tabular}	

\caption{Test 3: Comparison between the uncontrolled dynamics and the controlled one with $2$ and $3$ discrete controls \label{Table_comparison}}
						
	\end{table}

\subsection{Test 4: Dirichlet boundary control }\label{ssec:test4}
In the final example we deal with a boundary optimal control problem, in particular we consider a scalar control acting on the Dirichlet boundary condition. 
We choose the same initial condition fixed in Test 3, $T=1$, $\Delta t=0.1$ complemented with homogeneous Dirichlet boundary conditions for all the walls of the square, except for the top wall where we set the following condition
$$
u(t,x,\alpha)=g(x,t,\alpha), \quad (x_1,x_2)\in [0,1] \times \{1\}, t \in (0,T].
$$
To construct our reference trajectory, we run a simulation fixing $g(x,t,\alpha)=x(1-x)\sin{t}$ and we compute the corresponding numerical solution that we denote by $\{(\bm{\tilde{U}}^i,\bm{\tilde{V}}^i,\bm{\tilde{P}}^i)\}_{i=1}^{n_t}$.
Then, we set the optimal control problem considering the control set $A=[0,1]$ and  the following controlled boundary condition $g(x,t,\alpha)=x(1-x)\alpha(t)$.
In this case the aim of the control problem is to reach the final configuration of the pressure field $\bm{\tilde{P}}^{n_t}$, so we define the cost functional
$$
J(\alpha)= \| P^{n_t}(\alpha)-\bm{\tilde{P}}^{n_t}\|^2
$$
where $P^{n_t}(\alpha)$ is the solution of the optimal control problem at final time with control $\alpha$. Note that the running cost is 0 in this example.

The number of discrete controls in this example varies in the set $\{2,3,5\}$ and this will correspond to an increasing number of nodes in the tree. We want to examine the efficiency of the method in terms of the its pruning capacity and its accuracy in the approximation of the target solution.\\

The comparison of the performances of these three cases is reported in Table \ref{Table_comparison2}. Note that the cost functional is decreasing  to $O(10^{-7})$ as we increase the number of controls. In this example the pruning method is rather efficient, as we can notice by the $Ratio_p$ column. As we increase the parameter $M$, the value $Ratio_p$ gets one order of magnitude in each step.

	\begin{table}[htbp]							
		\centering
		\begin{tabular}{c |c c c}
					\toprule
				$M$	& Cost functional & Nodes & $Ratio_p$
					     \\
					\midrule
					$2$ & 1.00e-6 & 228 & 9 \\
					$3$  & 9.19e-7 & 710 & 125 \\
					$5$  &   1.74-7 & 4541 &  2.7e3\\
			\bottomrule		
		\end{tabular}	

\caption{Test 4: Comparison between the approximations of the optimal control problem varying the number of controls \label{Table_comparison2}}
						
	\end{table}
	
\begin{figure}[htb!]
\centering
\minipage{0.4\textwidth}
  \includegraphics[width=\linewidth]{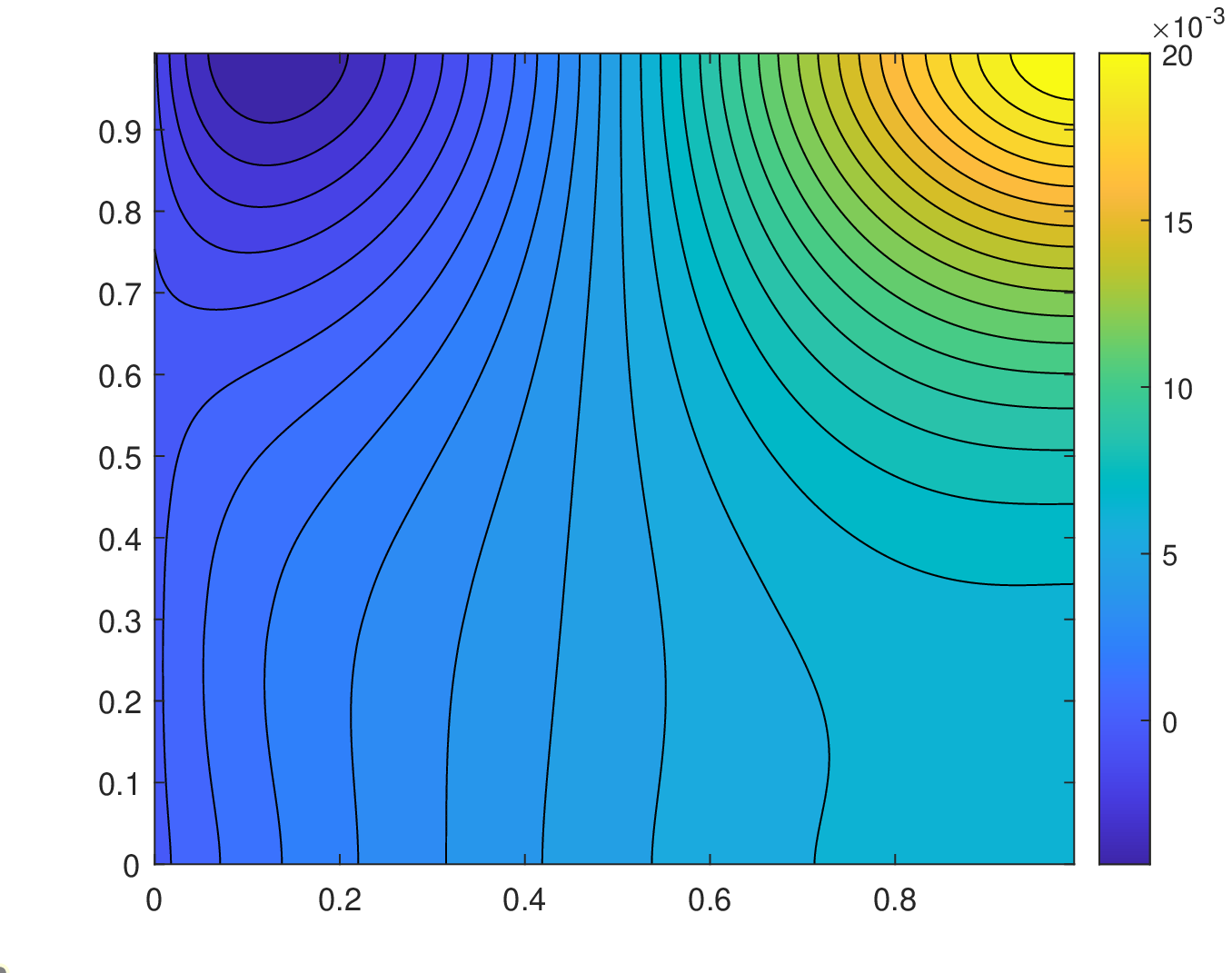}
\endminipage
\minipage{0.4\textwidth}
  \includegraphics[width=\linewidth]{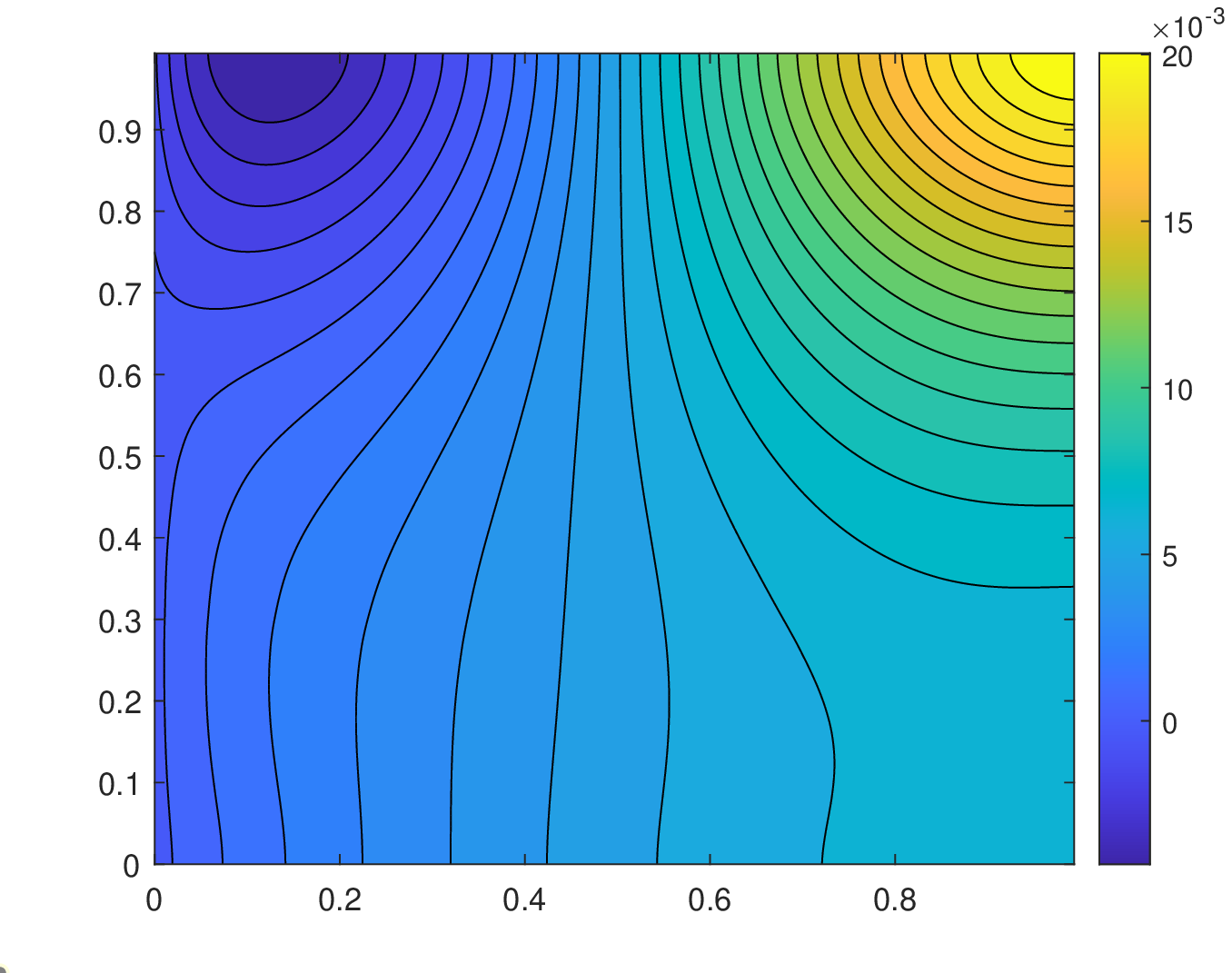}
\endminipage\hfill
\caption{Test 4: Final configuration of the reference solution (left) and final configuration of the controlled solution with $5$ discrete controls (right).}	
\label{fig_test_4}
\end{figure}

\begin{figure}[htb!]
\centering
\minipage{0.4\textwidth}
  \includegraphics[width=\linewidth]{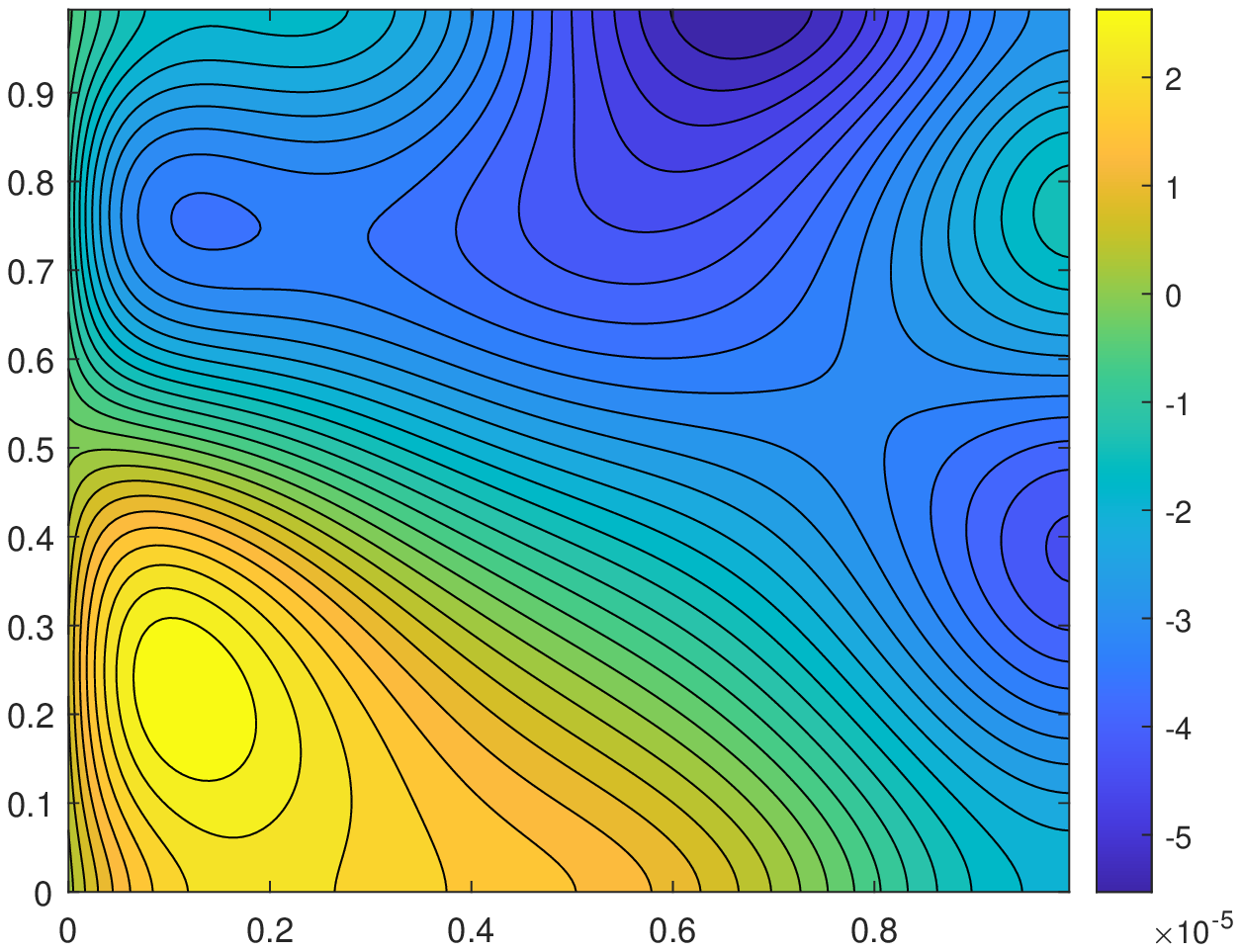}
\endminipage
\minipage{0.4\textwidth}
  \includegraphics[width=\linewidth]{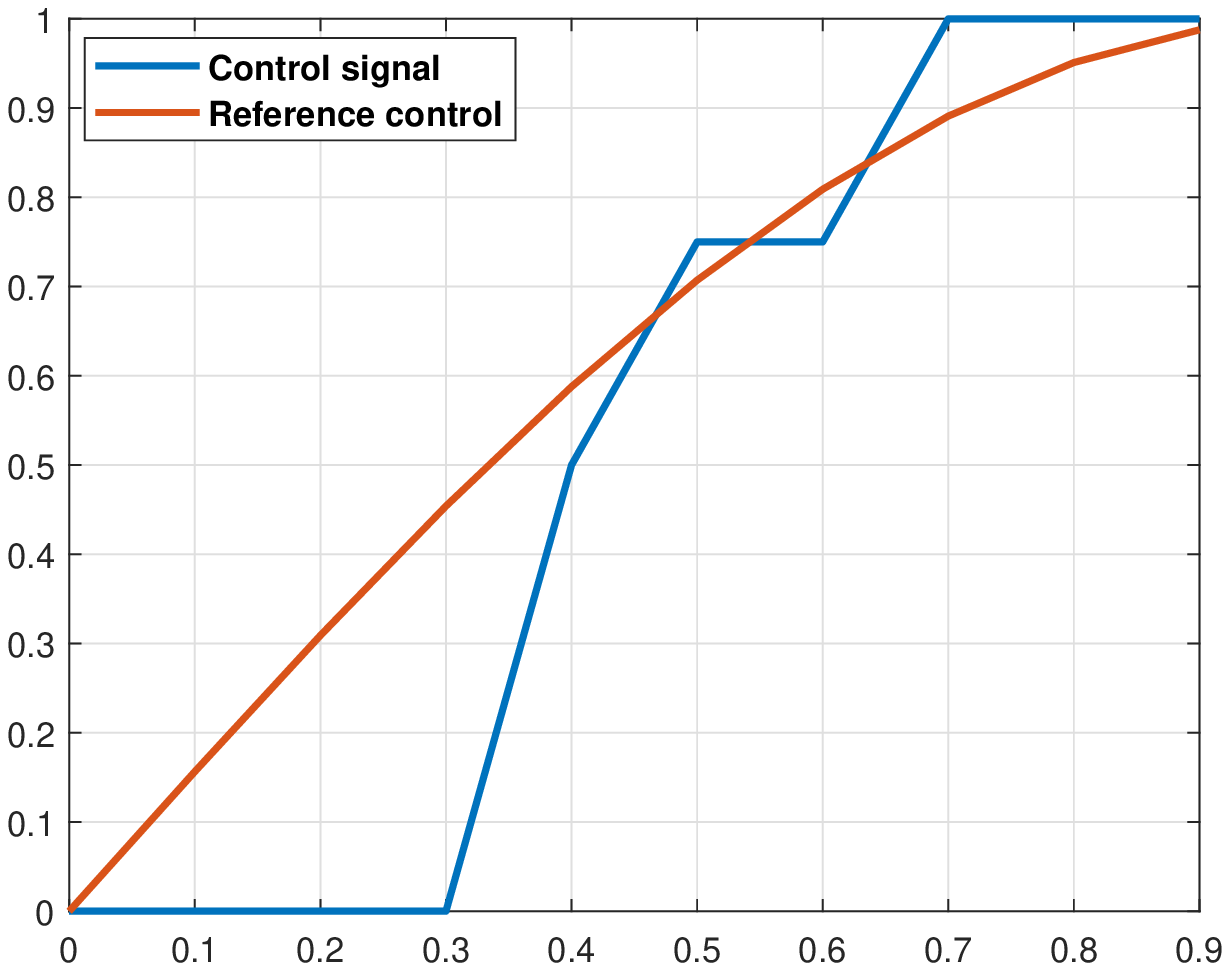}
\endminipage
\caption{Test 4: The difference between the reference solution and the controlled solution at final time (left) and comparison between the reference control and the computed optimal control (right).}	
\label{fig_test_4_2}
\end{figure}

In Figure \ref{fig_test_4} we show the configuration at final time of the reference solution $\bm{\tilde{P}}^{n_t}$ in the left panel and the controlled solution fixing $M=5$ of the controlled solution in the right panel. Visually they look very similar, but this is also certified in the left panel of Figure \ref{fig_test_4_2} showing that the difference between the reference and the controlled solution is order $O(10^{-5})$. This shows that the numerical method is able to reconstruct an optimal control driving the dynamics close to the reference solution. Finally, in the right panel of Figure \ref{fig_test_4_2} the reference control and the numerical approximation are shown, where we can note that the optimal control is trying to mimic the reference signal.

\section{Conclusions}
In this paper we have presented our first results on the numerical approximation of optimal control problems for the Navier-Stokes equation. The problem is discretized in space to obtain a system of ordinary differential equations, then we set the control problem on the finite dimensional system of ordinary differential equations corresponding to that semi-discretization. A crucial role is played by a very compact representation of the dynamical system and by a tree structure method to solve the problem via Dynamic Programming. More precisely, we have illustrated that by taking advantage of the rectangular domain, and a tensor-structured discretization basis, that the discrete NS equation can be written, integrated and reduced entirely in matrix form, to dramatically reduce the computational cost of integrating the discrete NS equation. On the other hand, the tree structure algorithm is used to counteract the curse of dimensionality arising from the optimal control problem and the Dynamic Programming.

The combination of these two methods shows that the Dynamic Programming approach can be used also in this area and that 
the synthesis  of optimal feedbacks can also be obtained for these huge optimization problems. This is a good omen for the future and we plan to investigate more in detail the convergence of feedback controls and other optimal control problems for fluids.

\printcredits
\vspace{0.2cm}\noindent
{\em Acknowledgments. The first author is a member of INDAM GNCS (Gruppo Nazionale di Calcolo Scientifico).
This research has been partially supported by the PRIN 2017 project "Innovative Numerical Methods for Evolutionary Partial Differential Equations and Applications", contract n. 2017KKJP4X.}

\bibliographystyle{cas-model2-names}

\bibliography{deimbib_3}

%

\end{document}